\documentclass[a4paper,12pt,reqno]{amsart}

\usepackage[headings]{fullpage}

\makeindex

\usepackage{amsfonts,graphics,amsmath,amsthm,amscd,amssymb,latexsym,euscript, enumerate}
\usepackage{epsfig,url}
\usepackage{flafter}
\usepackage[all,cmtip,line]{xy}
\usepackage{array}
\usepackage[english]{babel}
\usepackage{overpic}
\usepackage{subfig}
\usepackage{multirow}
\usepackage{microtype}
\usepackage[dvipsnames,svgnames,table]{xcolor}
\usepackage{tikz}
\usetikzlibrary{arrows.meta,positioning,shapes.geometric,calc}
\usepackage{tikz-cd}

\usepackage{wrapfig}
\usepackage[shortlabels]{enumitem}
\setlist[enumerate, 1]{$(1)$}
\usepackage[normalem]{ulem}

\usepackage{graphicx}
\usepackage{hyperref}
\hypersetup{
    colorlinks=true,    
    linkcolor=blue,          
    citecolor=blue,      
    filecolor=blue,      
    urlcolor=blue
}

\newtheorem{theorem}{Theorem}[section]
\newtheorem{lemma}[theorem]{Lemma}
\newtheorem{proposition}[theorem]{Proposition}
\newtheorem{conjecture}[theorem]{Conjecture}

\newtheorem*{theorem*}{Theorem}
\newtheorem*{example*}{Example}

\theoremstyle{plain}
\newtheorem{corollary}[theorem]{Corollary}

\theoremstyle{definition} 
\newtheorem{definition}[theorem]{Definition}
\newtheorem{definition-lemma}[theorem]{Definition-Lemma}

\newtheorem{example}[theorem]{Example}

\newtheorem{remark}[theorem]{Remark}

\numberwithin{equation}{section}

\newcommand{\R}{\mathbb{R}}
\newcommand{\Z}{\mathbb{Z}}

\newcommand{\Q}{\mathbb{Q}}
\newcommand{\mc}{\mathcal}
\def\P{\mathbb{P}}

\DeclareMathOperator{\ord}{ord}

\DeclareMathOperator{\Val}{Val}

\DeclareMathOperator{\Nef}{Nef}

\DeclareMathOperator{\Pic}{Pic}

\def\Pic{\operatorname{Pic}}

\def\mult{\operatorname{mult}}

\def\Nef{\operatorname{Nef}}

\def\Supp{\operatorname{Supp}}

\def\NE{\operatorname{NE}}

\def\lct{\operatorname{lct}}

\usepackage{mathtools}

\DeclarePairedDelimiterX{\norm}[1]{\lVert}{\rVert}{#1}

\title[Admissible extremal rays and the anticanonical minimal model program]
{Admissible extremal rays and the anticanonical minimal model program}

\begin{document}

\author[S. Choi]{Sung Rak Choi}
\author[S. Jang]{Sungwook Jang}
\author[D.-W. Lee]{Dae-Won Lee}
\address[Sung Rak Choi]{Department of Mathematics, Yonsei University, 50 Yonsei-ro, Seodaemun-gu, Seoul 03722, Republic of Korea}
\email{sungrakc@yonsei.ac.kr}
\address[Sungwook Jang]{Center for Complex Geometry, Institute for Basic Science, 34126 Daejeon, Republic of Korea}
\email{swjang@ibs.re.kr}
\address[Dae-Won Lee]{School of Mathematics, Korea Institute for Advanced Study, 85 Hoegiro, Dongdaemun-gu, Seoul 02455, Republic of Korea}
\email{daewonlee@kias.re.kr}


\date{\today}

\begin{abstract}
We develop an anticanonical minimal model program for projective $\Q$-factorial klt pairs $(X,\Delta)$ such that $-(K_X+\Delta)$ is pseudoeffective. We first prove the cone and contraction theorems for potentially klt triples and establish the existence of flips for such triples. Using these results, we introduce the notion of admissible extremal rays and show that $(K_X+\Delta)$-positive admissible extremal rays are precisely the rays that are to be contracted by $-(K_X+\Delta)$-MMP. The elementary MMP steps associated to such admissible rays preserve klt singularities. For potentially klt pairs $(X,\Delta)$ with $\Q$-divisors $\Delta$, we prove that every $(K_X+\Delta)$-positive extremal ray is admissible and obtain a cone theorem for the $(K_X+\Delta)$-positive part of the Mori cone. Assuming termination of the anticanonical minimal model program, we establish a dichotomy: potentially klt pairs admit anticanonical minimal models, whereas non-potentially klt pairs admit non-admissible models.
\end{abstract}

\maketitle


\section{Introduction}

Since the early stages of the development of the minimal model program, it has been recognized that working with pairs, rather than just varieties, is useful and essential for inductive purposes. Recent studies also show that the theory of generalized pairs has shown to be particularly powerful in birational geometry of algebraic varieties, especially in relative settings. Although they have been studied in various disguised forms for some time, the first rigorous definition was given by Birkar and Zhang in \cite{BZ}, inspired by the study of the effective Iitaka fibration conjecture. In this paper, we consider the \emph{triples} that are even more general than the generalized pairs. These triples were systematically studied in \cite{CJK}, \cite{CJKL}, where a triple was called a potential triple. By definition, a triple $(X,\Delta,D)$ consists of a pair $(X,\Delta)$ and a pseudoeffective $\R$-Cartier divisor $D$ on $X$. The singularities of triples, such as potentially klt and potentially lc, are defined by replacing the log discrepancy $a(E;X,\Delta)$ of a pair $(X,\Delta)$ by the log discrepancy $a(E;X,\Delta,D)\coloneqq a(E;X,\Delta)-\sigma_E(D)$ of a triple $(X,\Delta,D)$, where $\sigma_E(D)$ denotes the asymptotic divisorial valuation of $D$ along a prime divisor $E$ over $X$. For example, a triple $(X,\Delta,D)$ is  said to be potentially klt (or pklt for short) if $\inf_Ea(E;X,\Delta,D)>0$ holds where $\inf$ is taken over all prime divisors $E$ over $X$. See Section \ref{sect:prelim} for details.

We first prove the following cone theorem and contraction theorem which will provide a foundation for the minimal model program for pklt triples.

\bigskip

\begin{theorem}[Cone Theorem for pklt triples]\label{thrm:triple cone}
Let $(X,\Delta,D)$ be a pklt triple. Then the following hold:
\begin{enumerate}[$(1)$, leftmargin=*]
\item There exist at most countably many $(K_X+\Delta+D)$-negative extremal rays $R_i$ ($i\in I$) of the Mori cone $\overline{\NE}(X)$ such that
    $$
    \overline{\NE}(X)=\overline{\NE}(X)_{K_X+\Delta+D\geq0}+\sum_{i\in I}R_i.
    $$
    Furthermore, for any ample divisor $A$ on $X$, the set $I_A\coloneqq \{i\in I\mid (K_X+\Delta+D+A)\cdot R_i<0\}$ is finite and 
    $$
    \overline{\NE}(X)=\overline{\NE}(X)_{K_X+\Delta+D+A\geq0}+\sum_{i\in I_A}R_i.
    $$
\item Each $(K_X+\Delta+D)$-negative extremal ray $R$ is spanned by a rational curve $C$ such that $0<-(K_X+\Delta+D)\cdot C\leq 2\dim X$.
\end{enumerate}
\end{theorem}

The proof also shows that the family of $(K_X+\Delta+D)$-negative extremal rays of $\overline{\NE}(X)$ is locally finite 
at each $(K_X+\Delta+D)$-negative extremal ray $R_i$. Thus,   these extremal rays $R_i$ are discrete. These properties will be used to unveil the structure of the $(K_X+\Delta)$-positive half cone $\overline{\NE}(X)_{K_X+\Delta>0}$ for a klt pair $(X,\Delta)$ (Theorem \ref{thrm:pklt positive cone}).
Note that the classical cone theorem for pairs gives us no direct structural information on $\overline{\NE}(X)_{K_X+\Delta>0}$ or on the properties of $(K_X+\Delta)$-positive extremal rays.  The $(K_X+\Delta)$-positive part of $\overline{\NE}(X)$ can be circular by \cite[Theorem 2]{CRJ} or rational polyhedral (Example \ref{ex:only positive nonadmissible}).

\begin{theorem}[Contraction Theorem for pklt triples]\label{thrm:triple contraction}
Let $(X,\Delta,D)$ be a pklt triple, and let $R$ be a $(K_{X}+\Delta+D)$-negative extremal ray of $\overline{\NE}(X)$. Then there exists a projective contraction morphism $\varphi\colon X\to Y$ onto a normal projective variety such that $\varphi_{\ast}\mathcal{O}_X=\mathcal{O}_Y$, $\rho(X/Y)=1$ and an integral curve $C$ is contracted by $\varphi$ if and only if $[C]\in R$.
\end{theorem}

We remark that Theorem \ref{thrm:triple cone} and Theorem \ref{thrm:triple contraction} for the case $D=0$ recover the classical cone theorem and contraction theorem for klt pairs, respectively. The corresponding theorems for generalized pairs have also been established (\cite{HacLi}, \cite{Xie}).

One of the goals of studying triples is to establish the minimal model program for triples. As in the classical minimal model program, the following is an important ingredient for running the MMP on pklt triples.

\begin{theorem}\label{thrm:triple flip}
Let $(X,\Delta,D)$ be a $\Q$-factorial pklt triple. If the contraction $\varphi\colon X\to Y$ associated to a $(K_X+\Delta+D)$-negative extremal ray is of flipping type, then its $(K_X+\Delta+D)$-flip $\varphi^+\colon X^+\to Y$ exists. 
\end{theorem}

Combining the above results, we can run the minimal model program on pklt triples. The details are given in Section \ref{sect:triple MMP}, where we prove the following. 

\begin{theorem}[MMP on pklt triples]\label{thrm:MMP on pklt}
Let $(X,\Delta,D)$ be a $\Q$-factorial pklt triple. Then we can run the MMP on the triple $(X,\Delta,D)$.
\end{theorem}

The termination of the minimal model program on pklt triples naturally reduces to the termination of flips for pklt triples which still needs to be confirmed (Conjecture \ref{conj:pklt termination}). As the resulting models of this MMP, we expect to obtain minimal models or Mori fiber spaces of the given pklt triples. See Definition \ref{def:potential models} for the precise definitions of the models of triples. We note that many of the steps in the proofs of Theorems \ref{thrm:triple cone}, \ref{thrm:triple contraction}, and \ref{thrm:triple flip} require reductions to klt pairs. However, the termination of flips for pklt triples does not directly follow from the termination of flips for klt pairs. See Remark \ref{rem:insufficient pair termination}.

\medskip

As an important application of the MMP on pklt triples, we establish the $-(K_X+\Delta)$-MMP for a $\Q$-factorial klt pair $(X,\Delta)$ such that $-(K_X+\Delta)$ is pseudoeffective. For general klt pairs with pseudoeffective anticanonical divisor, the classical MMP does not directly provide such a program. The main reason is that many of the core ingredients that constitute the classical MMP are not applicable or they just fail to hold once we replace $K_X+\Delta$ by $-(K_X+\Delta)$. For example, the classical cone theorem  provides no information on the $-(K_X+\Delta)$-negative (that is, $(K_X+\Delta)$-positive) extremal rays of the Mori cone $\overline{\NE}(X)$. The fact that $-(K_X+\Delta)$-MMP worsens the singularities, possibly to non-lc, makes it even more difficult to run $-(K_X+\Delta)$-MMP (cf. \cite{GOST}, \cite{CP}). The anticanonical ring $R(-K_X)$ is not a birational invariant in general; preservation of this ring is nevertheless a natural requirement for the $-K_X$-MMP. These difficulties are often avoided by simply assuming that the varieties involved are Mori dream spaces. However, such an assumption is usually much stronger than needed since we are interested only in $-(K_X+\Delta)$-MMP. We note that there exists a variety $X$ that is not a Mori dream space and nevertheless admits a $-K_X$-MMP (Example \ref{ex:-K-MMP on nonMDS}).

We construct the $-(K_X+\Delta)$-MMP by first identifying the extremal rays of the Mori cone $\overline{\NE}(X)$ that are to be contracted. They are precisely the following type of extremal rays. For a given klt pair $(X,\Delta)$, an extremal ray $R$ of $\overline{\NE}(X)$ is called an \emph{admissible ray} for $(X,\Delta)$ if $R$ is $(K_X+\Delta+D)$-negative for some pklt triple $(X,\Delta,D)$ with a big $\R$-Cartier divisor $D$ (Definition \ref{def:admissible rays}). Assume in addition that $X$ is $\Q$-factorial and $-(K_X+\Delta)$ is pseudoeffective. The contraction morphism $\varphi_R\colon X\to Y$ associated to a $(K_X+\Delta)$-positive  admissible ray $R$ for $(X,\Delta)$ given by the Contraction Theorem (Theorem \ref{thrm:triple contraction}) induces a divisorial contraction or a flip.
We can regard such a birational map as an elementary birational map of $-(K_X+\Delta)$-MMP due to the following theorem.
We note that $\varphi_R$ cannot be a Mori fiber space since $-(K_X+\Delta)$ is pseudoeffective and the cone of movable curves is $(K_X+\Delta)$-nonpositive by \cite{BDPP}.

\begin{theorem}\label{thrm:elementary -K-MMP}
Let $(X,\Delta)$ be a $\Q$-factorial klt pair such that $-(K_X+\Delta)$ is pseudoeffective, and let the admissible extremal ray $R$ for $(X,\Delta)$, the triple $(X,\Delta,D)$ and the map $\varphi_R$ be as above.
\begin{enumerate}[leftmargin=*]
  \item If $\varphi_R$ is a divisorial contraction, then $(Y,\Delta_Y)$ is klt where $\Delta_Y\coloneqq \varphi_{R*}(\Delta)$ and $-(K_Y+\Delta_Y)$ is pseudoeffective.
  \item If $\varphi_R$ is a small contraction and $\psi\colon X\dashrightarrow X^+/Y$ is the flip of $\varphi_R$, then $(X^+,\Delta^+)$ is klt where $\Delta^+\coloneqq \psi_*(\Delta)$ and $-(K_{X^+}+\Delta^+)$ is pseudoeffective.
\end{enumerate}
\end{theorem}
\noindent These statements are proved in Theorem \ref{thrm:-KMMP preserv potdisc}.

\medskip

A $(K_X+\Delta)$-positive non-admissible ray $R$ for $(X,\Delta)$ can still be $(K_X+\Delta')$-negative for some klt pair $(X,\Delta')$ and such a ray $R$ can be contracted by a step of $(K_X+\Delta')$-MMP. However, we show in Proposition \ref{prop:non-admissible contraction} that such contractions are not admissible for the given pair $(X,\Delta)$ in the sense that the resulting pair of $(X,\Delta)$ is not klt. 
\medskip

\begin{figure}[!htbp]
\centering
\scalebox{0.92}{%
\begin{tikzpicture}[
  font=\linespread{1}\fontsize{10}{12}\selectfont,
  box/.style={draw,rounded corners=2pt,align=center,
    inner xsep=6pt,inner ysep=6pt,line width=.55pt},
  test/.style={draw,diamond,aspect=3,align=center,
    inner xsep=0pt,inner ysep=0pt,line width=.65pt},
  arrow/.style={-{Stealth[length=2.0mm,width=1.4mm]},line width=.55pt},
  edge label/.style={font=\linespread{1}\fontsize{9}{11}\selectfont}
]

\node[box,text width=9.3cm] (initial) at (0,0)
  {\textbf{START}: given a $\Q$-factorial klt pair $(X_1,\Delta_1)\coloneqq(X,\Delta)$\\
  such that $-(K_X+\Delta)$ is pseudoeffective.};
\node[test,text width=3.0cm] (nef) at (0,-2.0)
  {\boxed{\textbf{Q1}}\\Is $-(K_{X_i}+\Delta_i)$ nef?};
\node[box,text width=3.0cm] (minimal) at (6.6,-2.0)
  {\textbf{STOP}:\\ $(X_i,\Delta_i)$ is an anticanonical minimal model.};
\node[test,text width=4.4cm] (admissible) at (0,-4.9)
  {\boxed{\textbf{Q2}}\\
   {\fontsize{9}{11}\selectfont
    Is there a $(K_{X_i}+\Delta_i)$-positive\\
    admissible ray $R$?}};
\node[box,text width=3.0cm] (nonadmissible) at (6.6,-4.9)
  {\textbf{STOP}:\\$(X_i,\Delta_i)$ is a non-admissible model.};
\node[box,text width=6.8cm] (choose) at (0,-7.6)
  {There exists a birational contraction \\
   $\varphi_R\colon X_i\to Z_i$ associated to $R$\\
   by Theorem \ref{thrm:triple contraction}.};
\node[box,text width=4.3cm] (divisorial) at (-2.75,-10.4)
  {Let $X_{i+1}\coloneqq Z_i$\\
   $f_i\coloneqq \varphi_R$\\
   };
\node[box,text width=6.3cm] (flip) at (3.75,-10.4)
  {There exists a flip $f_i\colon X_i\dashrightarrow X_i^+ / Z_i$ \\
   by Theorem \ref{thrm:triple flip}.\\
   Let $X_{i+1}\coloneqq X_i^+$};
\node[box,text width=6.3cm] (update) at (0,-12.8)
  {Let $(X_{i+1},\Delta_{i+1})\coloneqq (X_{i+1}, f_{i*}\Delta_i)$};

\draw[arrow] (initial) -- (nef);

\draw[arrow] (nef.east)
  -- node[above,edge label]{Yes} (minimal.west);

\draw[arrow] (nef.south)
  -- node[right,edge label]{No} (admissible.north);

\draw[arrow] (admissible.east)
  -- node[above,edge label]{No} (nonadmissible.west);

\draw[arrow] (admissible.south)
  -- node[right,edge label]{Yes} (choose.north);

\coordinate (split) at (0,-9.3);

\draw (choose.south) -- (split);

\draw[arrow] (split) -|
  node[pos=.35,above,edge label]{$\varphi_R$ is divisorial}
  (divisorial.north);

\draw[arrow] (split) -|
  node[pos=.35,above,edge label]{$\varphi_R$ is small}
  (flip.north);

\coordinate (merge) at (0,-11.6);

\draw (divisorial.south) |- (merge);
\draw (flip.south) |- (merge);

\draw[arrow] (merge) -- (update.north);

\draw[arrow] (update.west) -- (-6.15,-12.8)
  -- node[left,edge label,rotate=90,anchor=south]{Return to \boxed{\textbf{Q1}}}
  (-6.15,-2.0) -- (nef.west);

\end{tikzpicture}%
}
\caption{The flow chart of $-(K_X+\Delta)$-MMP}
\label{fig:-K-MMP}
\end{figure}
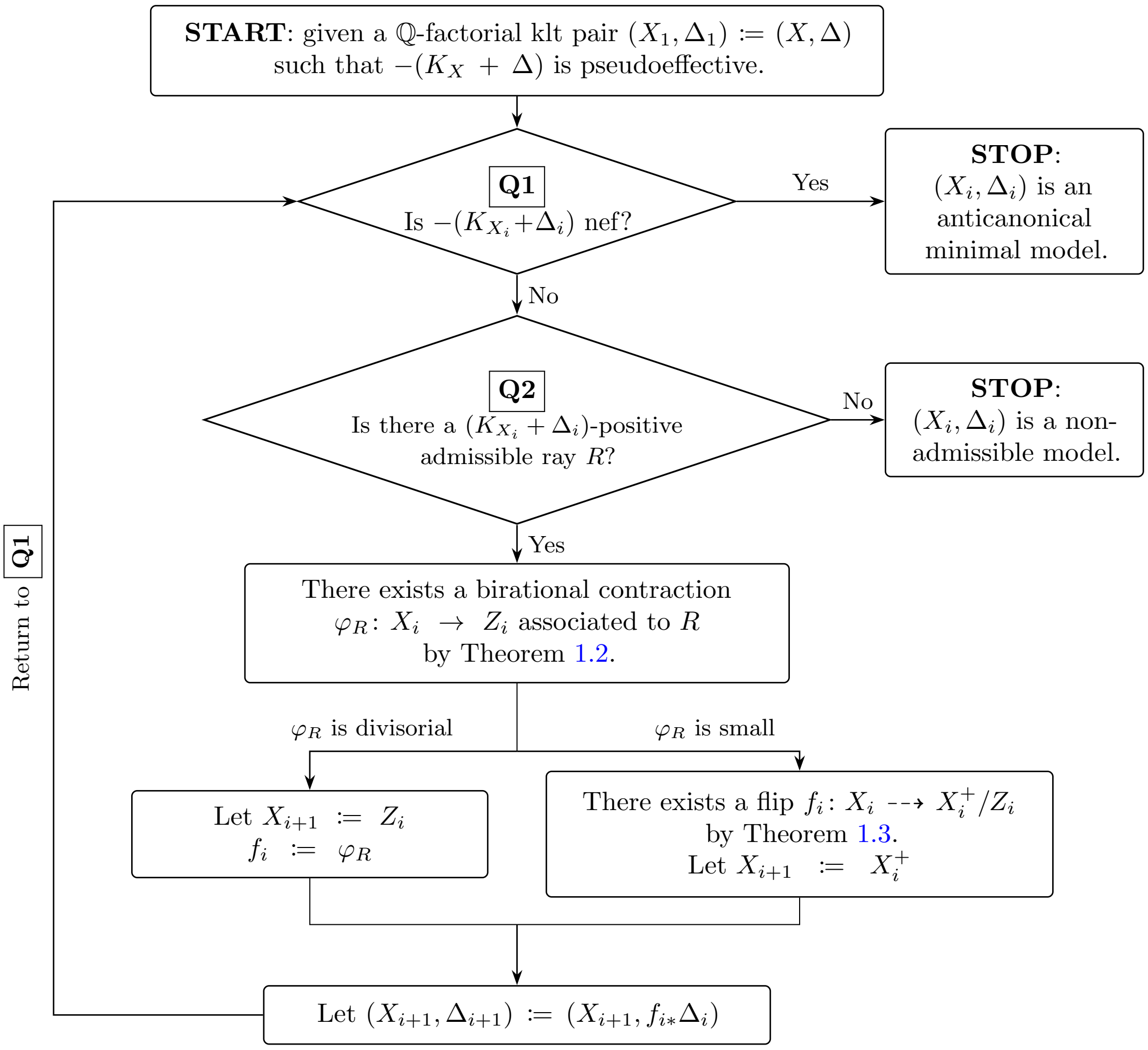

We run the $-(K_X+\Delta)$-MMP on a $\Q$-factorial klt pair $(X,\Delta)$ such that $-(K_X+\Delta)$ is pseudoeffective by successively applying the divisorial contractions and flips associated to $(K_X+\Delta)$-positive admissible rays.
The detailed flow chart of $-(K_X+\Delta)$-MMP is given in Figure \ref{fig:-K-MMP} of Section \ref{sect:klt -K MMP}. If the $-(K_X+\Delta)$-MMP terminates, then we obtain either an \textit{anticanonical minimal model} or a \emph{non-admissible model} (Definition \ref{def:nonadmissible model}).
See Example \ref{ex:only positive nonadmissible} for an example of a non-admissible model. In this example, $X$ is a Mori dream space having a unique $(K_X+\Delta)$-positive extremal ray which is non-admissible. This also shows that even if a $(K_X+\Delta)$-positive extremal ray $R$ is isolated and rational, it need not be admissible for $(X,\Delta)$.

As a by-product, a careful analysis of admissible rays yields a cone theorem for the $(K_X+\Delta)$-positive part of $\overline{\NE}(X)$ when the triple $(X,\Delta,D)$ with $D=-(K_X+\Delta)$ is pklt
(Theorem \ref{thrm:pklt positive cone}). This result together with the Cone Theorem for klt pairs implies that the
$(K_X+\Delta)$-nontrivial extremal rays can accumulate only toward the $(K_X+\Delta)$-trivial hyperplane in $N_1(X)_{\R}$
(Corollary \ref{cor:cone thrm for pklr pair}).

\medskip

This paper is organized as follows.
In Section \ref{sect:prelim}, we recall the basic definitions and results regarding the triples. Some basic properties of pklt triples are proved in this section. In Section 3, we establish the minimal model program on pklt triples by proving the cone theorem (Theorem \ref{thrm:triple cone}) and contraction theorem (Theorem \ref{thrm:triple contraction}) for pklt triples. It is constructed almost in parallel with the standard minimal model program. In Section \ref{sect:admissible rays}, we define and study the geometry of admissible rays. They are generalizations of $(K_X+\Delta)$-negative extremal rays for a klt pair $(X,\Delta)$. The structure theorem for the $(K_X+\Delta)$-positive part of the Mori cone $\overline{\NE}(X)$ is proved in terms of admissible rays. Lastly, in Section \ref{sect:klt -K MMP}, based on the properties of the admissible rays, we construct the $-(K_X+\Delta)$-minimal model program. The flowchart of $-(K_X+\Delta)$-minimal model program is given in Figure \ref{fig:-K-MMP} of this section.

\section*{Acknowledgement}
S. Choi and D.-W. Lee are partially supported by Samsung Science and Technology Foundation under Project Number SSTF-BA2302-03. D.-W. Lee is also partially supported by the National Research Foundation of Korea (No. RS-2023-00237440 and RS-2025-00513064) and by a KIAS Individual Grant (MG111201) at Korea Institute for Advanced Study.

\section{Preliminaries}\label{sect:prelim}

In this section, we collect definitions and basic results that will be used throughout the paper.

A variety is integral and is defined over an algebraically closed field of characteristic zero, e.g., the field of complex numbers $\mathbb C$.

A \emph{pair} $(X,\Delta)$ consists of a normal projective variety $X$ and an effective $\R$-divisor $\Delta$ on $X$ such that $K_X+\Delta$ is $\R$-Cartier.
For a birational morphism $f\colon Y\to X$, we can write $K_Y+\Delta_Y=f^{\ast}(K_X+\Delta)$ for some divisor $\Delta_Y$ on $Y$. For any prime divisor $E$ over $X$, by taking a higher model if necessary, we may assume that $E$ is a prime divisor on $Y$. Then the log discrepancy $a(E;X,\Delta)$ of $(X,\Delta)$ along $E$ is defined as $1-\mult_E \Delta_Y$. If the strict inequality $\inf_E a(E;X,\Delta)>0$ holds where $\inf$ is taken over all prime divisors $E$ over $X$, then the pair $(X,\Delta)$ is said to have klt singularities. If the inequality $\geq 0$ holds, then $(X,\Delta)$ is said to have lc singularities.

A \emph{triple} $(X,\Delta,D)$ consists of a pair $(X,\Delta)$ and a pseudoeffective $\R$-Cartier divisor $D$ on $X$.
This notion was previously studied in \cite{CJK}, \cite{CJKL}, \cite{CJL} where a triple was called a potential triple. For simplicity, we will call it a triple. The triples $(X,\Delta,D)$ with $D=-(K_X+\Delta)$ have been defined and studied in \cite{CP}.

\begin{definition}\label{def:div asymp valuation}
Let $D$ be a pseudoeffective $\R$-Cartier divisor on a normal projective variety $X$.
\begin{enumerate}[leftmargin=*]
\item  Let $E$ be a prime divisor over $X$. If $D$ is big and $f\colon Y\to X$ is a projective birational morphism with $Y$ smooth such that $E$ is a prime divisor on $Y$, then the \emph{asymptotic divisorial valuation $\sigma_E(D)$ of $D$ along $E$} is defined as 
$$\sigma_E(D)\coloneqq \inf \{\mult_E(D')\mid f^*D\equiv D'\geq 0\}.$$
For a pseudoeffective divisor $D$, we define $\sigma_E(D)\coloneqq \lim_{t\to 0^+} \sigma_E(D+tA)$, where $A$ is an ample divisor. This definition is independent of the choice of $f$ and $A$. 

\item Assume that $X$ is $\Q$-factorial. The \emph{negative part} $N_D$ of $D$ is defined as the divisor $N_D=\sum_E \sigma_E(D)E$ where $E$ runs over the prime divisors on $X$. (There are only finitely many prime divisors $E$ on $X$ such that $\sigma_E(D)>0$ by \cite{Nak}.) The \emph{positive part} $P_D$ of $D$ is defined as the numerical class $P_D\equiv D-N_D$. The decomposition $D\equiv P_D+N_D$ is called the \emph{divisorial Zariski decomposition of} $D$. If $P_D$ is nef, then it is called the \emph{Zariski decomposition}.
\end{enumerate}
\end{definition}

\begin{definition}\label{def:pklt triple}
Let $(X,\Delta,D)$ be a triple.
\begin{enumerate}[leftmargin=*]
\item The log discrepancy $a(E;X,\Delta,D)$ of $(X,\Delta,D)$ along a prime divisor $E$ over $X$ is defined as
$a(E;X,\Delta,D)\coloneqq a(E;X,\Delta)-\sigma_E(D)$
where $a(E;X,\Delta)$ is the log discrepancy of the pair $(X,\Delta)$ along $E$ and $\sigma_E(D)$ is the asymptotic divisorial valuation of $D$ along $E$.

\item A triple $(X,\Delta,D)$ is said to be \emph{potentially klt} (or \emph{pklt} for short) if $\inf_E a(E;X,\Delta,D)>0$ holds where $\inf$ is taken over all prime divisors $E$ over $X$. If $(X,\Delta,D)$ with $D=-(K_X+\Delta)$ is pklt, then we call $(X,\Delta)$ a \emph{pklt pair}.
\end{enumerate}
\end{definition}

Note that if $D$ is a nef divisor, then $\sigma_E(D)=0$ for any $E$ by \cite{ELMNP}. Thus, the triple $(X,\Delta,D)$ is pklt if and only if $(X,\Delta)$ is klt.

The following result is useful when we reduce the problems on pklt triples to the case of klt pairs.

\begin{theorem}[\protect{\cite[Theorem 1.1 and its proof]{CJK}}]\label{thrm:pklt to klt}
Let $(X,\Delta,D)$ be a pklt triple such that $D$ is a big $\Q$-divisor. Then there exists an effective divisor $D'\sim_\Q D$ such that $(X,\Delta+D')$ is a klt pair. Conversely, the existence of such an effective divisor $D'$ implies that the triple $(X,\Delta,D)$ is pklt.
\end{theorem}

For a triple $(X,\Delta,D)$ with $(X,\Delta)$ klt, the log canonical threshold $\lct_{\sigma}(X,\Delta,D)$ of $(X,\Delta,D)$ is defined as
$$
\lct_{\sigma}(X,\Delta,D)\coloneqq \inf_{v\in \Val_{X}^{\ast}}\frac{a(v;X,\Delta)}{\sigma_{v}(D)}=\inf_{E}\frac{a(E;X,\Delta)}{\sigma_{E}(D)}
$$
where $v$ is taken over all nontrivial valuations $v\in \Val_{X}^{\ast}$ on $X$ and $E$ is taken over all prime divisors $E$ over $X$ (cf. \cite{CJKL}). We define $\frac{a(v;X,\Delta)}{\sigma_{v}(D)}=+\infty$ if $\sigma_v(D)=0$, and $\frac{a(E;X,\Delta)}{\sigma_{E}(D)}=+\infty$ if $\sigma_E(D)=0$.

\begin{proposition}[\protect{\cite[Proposition 4.1]{CJKL}, \cite[Lemma 2.5]{CJK}}]\label{prop:pklt as lct>1}
Let $(X,\Delta)$ be a klt pair with a $\Q$-divisor $\Delta$ and let $D$ be a pseudoeffective $\Q$-Cartier $\Q$-divisor. Then $(X,\Delta,D)$ is pklt if and only if $\lct_{\sigma}(X,\Delta,D)>1$. 
\end{proposition}

\noindent For a $\Q$-factorial klt pair $(X,\Delta)$ and a pseudoeffective divisor $D$, we also consider the following pklt threshold:
$$\lambda_D\coloneqq \sup\{\lambda\geq 0\mid (X,\Delta,\lambda D)\;\text{is pklt}\}.$$

\begin{lemma}\label{lem:pklt<plct}
Let $(X,\Delta)$ be a $\Q$-factorial klt pair and $D$ be a pseudoeffective $\R$-Cartier divisor on $X$.
Then $(X,\Delta,sD)$ is pklt for any real number $s$ such that $0\leq s <\lct_{\sigma}(X,\Delta,D)$. In particular, we have $\lambda_D=\lct_{\sigma}(X,\Delta,D)$.
\end{lemma}
\begin{proof}
First of all, let $\delta\coloneqq \inf\{ a(E;X,\Delta)\mid \text{$E$ is a prime divisor over $X$}\}>0$. Let $s$ be a real number such that $0\leq s<\lct_{\sigma}(X,\Delta,D)$. Then for some fixed $t$ satisfying $s<t<\lct_{\sigma}(X,\Delta,D)$, we have $a(E;X,\Delta,tD)\geq 0$ for any prime divisor $E$ over $X$. We note that the following inequality holds for any prime divisor $E$ over $X$:
\begin{align*}
a(E;X,\Delta,sD)&=\left(1-\frac{s}{t}\right)a(E;X,\Delta)+\frac{s}{t}(a(E;X,\Delta)-t\sigma_E(D))\\
&\geq \left(1-\frac{s}{t}\right)\delta>0.
\end{align*}
This shows that $(X,\Delta,sD)$ is pklt. This also implies that $\lambda_D\geq \lct_{\sigma}(X,\Delta,D)$. The reverse inequality $\lambda_D\leq \lct_{\sigma}(X,\Delta,D)$ follows by definition.
\end{proof}

\begin{theorem}\label{thrm:pklt open}
Let $(X,\Delta,D)$ be a $\Q$-factorial pklt triple with $\Q$-divisors $\Delta$ and $D$.
\begin{enumerate}[$(1)$, leftmargin=*]
\item If the triples $(X,\Delta,D_1)$ and $(X,\Delta,D_2)$ with pseudoeffective $\R$-divisors $D_i \;(i=1,2)$ are pklt, then $(X,\Delta,\alpha D_1+\beta D_2)$ is pklt for any $\alpha,\beta\geq 0$ such that $\alpha+\beta=1$.

\item  For a sufficiently small $\varepsilon>0$, $(X,\Delta,(1+\varepsilon)D)$ is pklt.

\item If $D$ is big, then there is a neighborhood $U$ of $[D]$ in $N^1(X)_{\R}$ such that $(X,\Delta,D')$ is pklt for any divisor $D'$ with $[D']\in U$. 

\item  For every effective divisor $\Delta'$, the triple $(X,\Delta+\varepsilon\Delta',D)$ is pklt for all sufficiently small $\varepsilon\geq 0$.
\end{enumerate}
\end{theorem}

\begin{proof}
(1) Note that the inequality $\alpha\sigma_E(D_1)+\beta \sigma_E(D_2)\geq \sigma_E(\alpha D_1+\beta D_2)$ holds for any prime divisor $E$ over $X$ and for any $\alpha,\beta\geq 0$ such that $\alpha+\beta=1$.
Thus, we obtain
$$
a(E;X,\Delta,\alpha D_1+\beta D_2)\geq \alpha a(E;X,\Delta,D_1) + \beta a(E;X,\Delta,D_2).
$$
Taking the infimum of both sides over the prime divisors $E$ over $X$ proves the statement.

\noindent (2) This is \cite[Proposition 4.1]{CJKL}.

\noindent (3) We note first that $\gamma(B)\coloneqq \sup_E \frac{\sigma_E(B)}{a(E;X,\Delta)}<\infty$ for any big $\Q$-Cartier divisor $B$ by \cite[Proof of Theorem 2.4]{KL1}. We can also easily check that $(X,\Delta,B)$ is pklt if and only if $\gamma(B)<1$. Hence, for all sufficiently small $t> 0$, we have
\begin{equation}\label{eqn:KL}
\gamma(D+tB)\leq \gamma(D)+t\gamma(B)<1
\end{equation}
and this implies that $(X,\Delta,D+tB)$ is pklt.

Let $r\coloneqq \rho(X)\geq 1$.  
We can take a basis $\{[A_1],\ldots,[A_{r}]\}$ consisting of the classes of big $\Q$-Cartier $\Q$-divisors $A_1,\ldots,A_{r}$ such that $[D]=\sum_{i=1}^r a_i [A_i]$ for some $a_i>0$ for all $i$. By taking $A_i$ sufficiently small, we may also assume that $D+A_i$ are big and $(X,\Delta,D+A_i)$ are pklt by (\ref{eqn:KL}).
It is clear that $(X,\Delta,cD)$ is pklt for any $c$ such that $0<c<1$. Then for some fixed such $c$, the convex hull $\mathcal{C}$ in $N^1(X)_{\R}$ of the classes $[cD],[D+A_1],\cdots,[D+A_r]$ is $r$-dimensional and $(X,\Delta,D')$ is pklt whenever $[D']\in \mathcal{C}$ by (1). Furthermore, $[D]$ belongs to the interior of $\mathcal{C}$. Therefore, any open subset $U\subseteq \mathcal{C}$ such that $[D]\in U$ satisfies the condition in the statement.

\noindent (4) Let $\Delta'$ be an effective divisor on $X$. Since $(X,\Delta)$ is klt, there is a positive real number $\mu>0$ such that $(X,\Delta+\mu\Delta')$ is klt, i.e., $a(E;X,\Delta+\mu\Delta')>0$ for any prime divisor $E$ over $X$. Thus, for a sufficiently small real number $\varepsilon'>0$, we have
\begin{align*}
    a(E;X,\Delta+\mu\varepsilon'\Delta',D)&=a(E;X,\Delta+\mu\varepsilon'\Delta')-\sigma_{E}(D)\\
    &=a(E;X,(1-\varepsilon')\Delta+\varepsilon'(\Delta+\mu\Delta'))-\sigma_E(D)\\
    &= (1-\varepsilon')a(E;X,\Delta)+\varepsilon'a(E;X,\Delta+\mu\Delta')-\sigma_E(D)\\
    &> (1-\varepsilon')a(E;X,\Delta)-\sigma_{E}(D)\\
    &=(1-\varepsilon')a(E;X,\Delta,\tfrac{1}{1-\varepsilon'}D).
\end{align*}
Since we may assume that $(X,\Delta,\tfrac{1}{1-\varepsilon'}D)$ is pklt by (2), there exists a constant $r>0$ such that $a(E;X,\Delta,\tfrac{1}{1-\varepsilon'}D)>r$ for any prime divisor $E$ over $X$. Therefore, we obtain the inequality
$$a(E;X,\Delta+\mu\varepsilon'\Delta',D)\geq (1-\varepsilon')r>0,$$
which implies that $(X,\Delta+\mu\varepsilon'\Delta',D)$ is pklt. This proves the claim since we can take the real number $\mu\varepsilon'$ arbitrarily small.
\end{proof}

\medskip

Note that in Theorem \ref{thrm:pklt open}, the divisors $\Delta,D$ in $(X,\Delta,D)$  are $\Q$-divisors since the results and arguments used in the proof are valid for $\Q$-divisors. However, it is natural to expect that all the statements hold for $\R$-divisors as well. See \cite{CJK}, \cite{CJKL} for detailed expositions on the triples.

A birational map $f\colon X\dashrightarrow Y$ is a \emph{contraction} if the inverse $f^{-1}$ does not contract divisors. Let $M$ and $M_Y\coloneqq f_*M$ be $\R$-Cartier divisors on $X$ and $Y$, respectively. We say that $f$ is $M$-\emph{nonpositive} if, on a common resolution $p\colon W\to X$ and $q\colon W\to Y$, we have
$$
p^{\ast}M=q^{\ast}M_Y+F
$$
for an effective $q$-exceptional divisor $F$. We say $f$ is $M$-\emph{negative} if $\Supp F$ contains the strict transforms of all prime divisors on $X$ contracted by $f$.

\begin{definition}\label{def:potential models}
Let $(X,\Delta,D)$ be a pklt triple and let $f\colon X\dashrightarrow Y$ be a birational contraction, with $Y$ $\Q$-factorial. Let $\Delta_Y\coloneqq f_{\ast}\Delta$ and $D_Y\coloneqq f_{\ast}D$, and let $(Y,\Delta_Y,D_Y)$ be a pklt triple.

\begin{enumerate}[$(1)$, leftmargin=*]
\item The triple $(Y,\Delta_Y,D_Y)$ is a \emph{minimal model} of $(X,\Delta,D)$ if $K_Y+\Delta_Y+D_Y$ is nef and the inequality
\begin{equation*}
    a(E;X,\Delta,D)\leq a(E;Y,\Delta_Y,D_Y)
\end{equation*}
holds for any prime divisor $E$ over $X$, with the strict inequality for every prime divisor $E$ on $X$ contracted by $f$.
\item The map $f$ with a contraction $g\colon Y\to Z$ defines a \emph{Mori fiber space} of $(X,\Delta,D)$ if $\dim Z<\dim Y$, $\rho(Y/Z)=1$, $-(K_Y+\Delta_Y+D_Y)$ is relatively ample over $Z$, and the inequality
\begin{equation*}
    a(E;X,\Delta,D)\leq a(E;Y,\Delta_Y,D_Y)
\end{equation*}
holds for any prime divisor $E$ over $X$, with strict inequality for every prime divisor $E$ on $X$ contracted by $f$.
\end{enumerate}
\end{definition}

Note that if $(X,\Delta,D)$ is a $\Q$-factorial pklt triple with $D=-(K_X+\Delta)$, then it is trivially a minimal model.

As in the case of pairs and generalized pairs, we expect that the above models are obtained as resulting models of the minimal model program on triples, which will be developed in Section \ref{sect:triple MMP}.

\begin{definition}\label{def:triple flip}
For $(X,\Delta,D)$ a $\Q$-factorial triple, an \emph{$(X,\Delta,D)$-flipping contraction} is a small projective contraction $g\colon X\to Z$ such that $\rho(X/Z)=1$ and $-(K_X+\Delta+D)$ is $g$-ample. 
The \emph{flip} of $g$ is a small projective contraction $g^{+}\colon X^{+}\to Z$ with the induced isomorphism in codimension one $f\colon X\dashrightarrow X^+$ such that $X^{+}$ is $\Q$-factorial, $\rho(X^{+}/Z)=1$, and $K_{X^{+}}+\Delta^{+}+D^{+}$ is $g^+$-ample where $\Delta^{+}\coloneqq f_{\ast}\Delta$ and $D^{+}\coloneqq f_{\ast}D$.
\end{definition}

\begin{lemma}[\protect{cf. \cite[Theorem 3.7]{KM}}]\label{lem:equiv to sim}
Let $(X,\Delta)$ be a $\Q$-factorial klt pair and let $g\colon X\to Z$ contract a $(K_X+\Delta)$-negative extremal ray $R$. 
\begin{enumerate}[$(1)$, leftmargin=*]
\item If an $\R$-Cartier divisor $M$ is numerically trivial on $R$, then $M\sim_{\R} g^{\ast}M_Z$ for an $\R$-Cartier divisor $M_Z$ on $Z$. 

\item If $g$ is small and $g^{+}\colon X^{+}\to Z$ is its flip, then for any $\R$-Cartier divisor $M$ on $X$ with $M\not\equiv_Z0$, $-M$ is $g$-ample if and only if its strict transform $M^+$ is $g^+$-ample.
\end{enumerate}
\end{lemma}
\begin{proof}
(1) The $\Q$-Cartier divisor case follows from \cite[Theorem 3.7]{KM}. Let $C$ be a curve spanning the ray $R$. By \cite[Theorem 3.7]{KM} (see also \cite[Corollary 3.17]{KM}), we have an exact sequence
$$
\begin{tikzcd}
\Pic(Z) \arrow[r,"L\mapsto g^{\ast}L"] & \Pic(X) \arrow[r, "M\mapsto M\cdot C"] & \Z.
\end{tikzcd}
$$
Since $\R$ is flat over $\Z$, by tensoring the above sequence with $\R$, we obtain an exact sequence
$$
\begin{tikzcd}
\Pic(Z)_{\R} \arrow[r] & \Pic(X)_{\R} \arrow[r] & \R.
\end{tikzcd}
$$
This implies the statement (1).

\noindent (2) Suppose that $g$ is small and $g^+\colon X^+\to Z$ is its flip. Fix again an integral curve $C$ spanning the ray $R$. Then since $\rho(X/Z)=1$, there exists a real number $c\neq 0$ such that $M- c(K_X+\Delta)\equiv 0$ over $Z$. By (1), there exists an $\R$-Cartier divisor $T$ on $Z$ such that $M-c(K_X+\Delta)\sim_\R g^*T$. Taking strict transforms under $X\dashrightarrow X^+$, we obtain
$$
M^+-c(K_{X^+}+\Delta^+)\sim_\R (g^{+})^*T
$$
where $\Delta^+$ is the strict transform of $\Delta$ on $X^+$.
By the definition of the flip, $-(K_X+\Delta)$ is $g$-ample and $(K_{X^+}+\Delta^+)$ is $g^+$-ample. Thus, the following three conditions are equivalent: (i) $-M$ is $g$-ample, (ii) $c>0$, and (iii) $M^+$ is $g^+$-ample.
\end{proof}
The existence of flips of a general pklt triple is proved in Theorem \ref{thrm:pklt flip exists}.

\begin{proposition}[{cf. \cite[Lemma 2.5]{CJL}}]\label{prop:pdisc ineq}
Let $(X,\Delta,D)$ be a $\Q$-factorial pklt triple. Assume that the birational contraction $f\colon X\dashrightarrow Y$ to a $\Q$-factorial variety $Y$ is a divisorial contraction or a flip of $(X,\Delta,D)$ associated to a $(K_X+\Delta+D)$-negative extremal ray $R$. Let $\Delta_Y\coloneqq f_{\ast}\Delta$ and $D_Y\coloneqq f_{\ast}D$.
Then for every prime divisor $E$ over $X$
$$
a(E;X,\Delta,D)\leq a(E;Y,\Delta_Y,D_Y)
$$
holds. In particular, $(Y,\Delta_Y,D_Y)$ is a $\Q$-factorial pklt triple. Moreover, the inequality is strict for every $f$-exceptional prime divisor $E$ on $X$.
\end{proposition}
\begin{proof}
The $\Q$-factoriality is preserved since each divisorial contraction or flip is constructed as a step of an MMP for a suitable klt pair.

Let $p\colon W\to X$ and $q\colon W\to Y$ be a common resolution on which $E$ appears as a prime divisor. 
Since $f$ is $(K_X+\Delta+D)$-negative, the negativity lemma gives
\begin{equation}\label{eqn:**}\tag{$**$}
p^*(K_X+\Delta+D)=q^*(K_{Y}+\Delta_Y+D_Y)+F
\end{equation}
for some effective divisor $F$ on $W$.

If $D\cdot R\leq 0$, then we may write $p^*D=q^*D_Y+G$ for some effective $q$-exceptional divisor $G$.
Subtracting this equality from (\ref{eqn:**}) gives us 
$$
p^*(K_X+\Delta) = q^*(K_{Y}+\Delta_Y)+F-G.
$$
Hence, we have $a(E;Y,\Delta_Y)-a(E;X,\Delta) = \mult_{E}(F-G)$.  Using the equality $\sigma_E(D)=\sigma_E(D_Y)+\mult_{E} G$, we also obtain
$$
\begin{aligned}
a(E;Y,\Delta_Y,D_Y)-a(E;X,\Delta,D)
&= a(E;Y,\Delta_Y)-a(E;X,\Delta)
+\sigma_E(D)-\sigma_E(D_Y) \\
&= \mult_{E}(F-G)+\operatorname{mult}_E G \\
&= \mult_{E} F.
\end{aligned}
$$

If $D\cdot R>0$, then $q^*D_Y=p^*D+G$ for some effective divisor $G$. Combining this with (\ref{eqn:**}), we obtain
$$
p^*(K_X+\Delta)=q^*(K_Y+\Delta_Y)+F+G.
$$
Thus, $a(E;X,\Delta)=a(E;Y,\Delta_Y)-\mult_E(F+G)$ holds for any prime divisor $E$ over $X$. Since $\sigma_E(D_Y)\leq \sigma_E(D)+\mult_EG$, we have
$$
a(E;Y,\Delta_Y,D_Y)-a(E;X,\Delta,D)\geq\mult_EF.
$$

Note that regardless of whether $R$ is $D$-nonpositive or $D$-positive, we have
$$
a(E;Y,\Delta_Y,D_Y)-a(E;X,\Delta,D)\geq\mult_EF.
$$
Since $\Supp F$ contains all the strict transforms of all $f$-exceptional prime divisors on $X$, the last statement follows.
\end{proof}

\medskip

\section{Minimal Model Program on potential triples}\label{sect:triple MMP}

In this section, we prove our first main theorems of the paper. Based on these theorems, we construct the MMP on pklt triples.

\begin{proof}[Proof of Theorem \ref{thrm:triple cone} \textrm{(Cone Theorem for pklt triples)}]
(1) Let $(X,\Delta,D)$ be a pklt triple with a pseudoeffective divisor $D$.
For any ample divisor $A$ on $X$, we have
\begin{align*}
\overline{\NE}(X)&=\overline{\NE}(X)_{K_X+\Delta+D+A\geq0}+\overline{\sum_{i\in I_A} R_i}
\end{align*}
where $R_i$ ($i\in I_A$) are $(K_X+\Delta+D+A)$-negative extremal rays of $\overline{\NE}(X)$. We can decompose $A=H+H'$ into two ample divisors $H,H'$ such that $D+H$ is a big $\Q$-divisor.
Since we may assume that $(X,\Delta,D+H)$ is a pklt triple, by Theorem \ref{thrm:pklt to klt} there exists an effective divisor $D'\sim_{\R}D+H$ such that $(X,\Delta+D')$ is a klt pair. Since $K_X+\Delta+D'+H'\equiv K_X+\Delta+D+A$, the set $\{R_i\mid i\in I_A\}$ of $(K_X+\Delta+D+A)$-negative extremal rays is finite by the Cone Theorem for klt pairs. This implies the following:
\begin{align*}
\overline{\NE}(X)&=\overline{\NE}(X)_{K_X+\Delta+D'+H'\geq 0}+\sum_{i\in I_A} R_i\\
&=\overline{\NE}(X)_{K_X+\Delta+D+A\geq0}+\sum_{i\in I_A} R_i.
\end{align*}
Note that an extremal ray $R$ is $(K_X+\Delta+D)$-negative  if and only if it is $(K_X+\Delta+D+qA)$-negative for some small positive rational number $q$. Since each set $\{R_i\mid i\in I_{qA}\}$ is finite, there exist at most countably many $(K_X+\Delta+D)$-negative extremal rays.

(2) As we just argued in the proof of (1), any $(K_X+\Delta+D)$-negative extremal ray $R$ is $(K_X+\Delta')$-negative for some klt pair $(X,\Delta')$.
Therefore, $R$ is spanned by a rational curve $C$ by the cone theorem for klt pairs. Now fix an ample Cartier divisor $H$ and fix a rational curve $C$ spanning the ray $R$ of minimal degree with respect to $H$. Such a curve exists since these degrees are positive integers. Let $H_i (i\geq 1)$  be a sequence of ample divisors such that $H_i\to 0$ and $D+H_i$ is a $\Q$-divisor for each $i$. By taking sufficiently small ample divisors $H_i$, we may assume that the ray $R$ is $(K_X+\Delta+D+H_i)$-negative. We may also assume that $(X,\Delta,D+H_i)$ is pklt and thus, there exists an effective divisor $G_i\sim_{\R} D+H_i$ such that $(X,\Delta+G_i)$ is klt. The Cone Theorem for the klt pair $(X,\Delta+G_i)$ gives a rational curve $C_i$ spanning $R$ such that
$$
0<-(K_X+\Delta+D+H_i)\cdot C_i\leq2\dim X.
$$
Write $[C_i]=u_i[C]$ in $N_1(X)_{\R}$ where $u_i\coloneqq \frac{H\cdot C_i}{H\cdot C}\geq1$. Then
$$
0<-(K_X+\Delta+D+H_i)\cdot C
=\frac{-(K_X+\Delta+D+H_i)\cdot C_i}{u_i}
\leq2\dim X.
$$
The curve $C$ is independent of $i$. Taking the limit $i\to \infty$ gives the desired inequality $0<-(K_X+\Delta+D)\cdot C\leq2\dim X$, which proves (2).
\end{proof}
\medskip

Now we prove the contraction theorem for triples.

\begin{proof}[Proof of Theorem \ref{thrm:triple contraction}\textrm{ (Contraction Theorem for pklt triples)}] 
Choose a sufficiently small ample $\R$-divisor $H$ on $X$ such that $D+H$ is a $\Q$-divisor and $(K_X+\Delta+D+H)\cdot R<0$. Since we may assume that the triple $(X,\Delta,D+H)$ is pklt, there exists an effective divisor $G$ such that $G\sim_{\Q} D+H$ and $(X,\Delta+G)$ is klt by Theorem \ref{thrm:pklt to klt}. Then $R$ is $(K_X+\Delta+G)$-negative. The contraction theorem for klt pairs \cite[Theorem 3.7]{KM} gives the required contraction.
\end{proof}

\begin{theorem}\label{thrm:pklt flip exists}
Let $(X,\Delta,D)$ be a $\Q$-factorial pklt triple and $R$ be a $(K_X+\Delta+D)$-negative extremal ray of $\overline{\NE}(X)$. Suppose that the contraction $\varphi_R=\varphi\colon X\to Y$ associated to $R$ is a $(X,\Delta,D)$-flipping contraction. Then the flip of $\varphi$ exists.
\end{theorem}
\begin{proof}
Let $A$ be a small ample divisor on $X$ such that $R$ is $(K_X+\Delta+D+A)$-negative. We may assume that $(X,\Delta,D+A)$ is pklt and $D+A$ is a big $\Q$-divisor.
By Theorem \ref{thrm:pklt to klt}, there exists an effective divisor $D'\sim_\Q D+A$ such that $(X,\Delta+D')$ is a klt pair. Therefore, $\varphi\colon X\to Y$ can be considered as a  $(K_X+\Delta+D')$-negative extremal small contraction associated to $R$. By \cite{bchm}, there exists a flip for $\varphi$ of the klt pair $(X,\Delta+D')$:
\begin{equation}\label{eqn:flip}
\begin{tikzcd}
X \arrow[rr,dashed,"f"] \arrow[dr,"\varphi"']
&& X^+ \arrow[dl,"\varphi^+"] \\
& Y. &
\end{tikzcd}
\end{equation}

We claim that this is also the flip for the pklt triple $(X,\Delta,D)$ associated to $R$.
Note that since $\rho(X/Y)=1$, there exists a positive real number $r>0$ such that 
$$
K_X+\Delta+D\equiv r(K_X+\Delta+D+A)\equiv r(K_X+\Delta+D') \;\text{over }Y.
$$
By Lemma \ref{lem:equiv to sim}, there exists an $\R$-Cartier divisor $T$ on $Y$ such that
$$
K_X+\Delta+D\sim_\R r(K_X+\Delta+D')+\varphi^*T.
$$
Taking strict transforms under $f$ in \ref{eqn:flip} gives 
$$
f_*(K_X+\Delta+D)\sim_\R rf_*(K_X+\Delta+D')+(\varphi^{+})^*T. 
$$
where $f_*(K_X+\Delta+D')$ is ample over $Y$. Hence, $f_*(K_X+\Delta+D)$ is ample over $Y$ and this shows that $X\dashrightarrow X^+ /Y$ (\ref{eqn:flip}) is also the flip for the pklt triple $(X,\Delta,D)$.
\end{proof}

\medskip

Now applying the cone theorem (Theorem \ref{thrm:triple cone}) and contraction theorem (Theorem \ref{thrm:triple contraction}) for pklt triples together with the existence of pklt flips (Theorem \ref{thrm:pklt flip exists}), we run the MMP on a $\Q$-factorial pklt triple as follows.

Let $(X,\Delta,D)=(X_1,\Delta_1,D_1)$ be a $
\Q$-factorial pklt triple and suppose that we have constructed $(X_i,\Delta_i,D_i)$ $(i\geq 1)$. If $K_{X_i}+\Delta_i+D_i$ is nef, then $(X_i,\Delta_i,D_i)$ is a minimal model of $(X,\Delta,D)$. Suppose that $K_{X_i}+\Delta_i+D_i$ is not nef. Then by the cone theorem (Theorem \ref{thrm:triple cone}), there exists a $(K_{X_i}+\Delta_i+D_i)$-negative extremal ray $R$ of $\overline{\NE}(X_i)$. By the contraction theorem (Theorem \ref{thrm:triple contraction}), there exists a contraction $\varphi_R\colon X_i\to X'_i$ associated to $R$. 

\begin{enumerate}
\item If $\varphi_R$ is divisorial (i.e., the $\varphi_R$-exceptional locus has codimension one), then by Proposition \ref{prop:pdisc ineq}, $(X'_i,\varphi_{R*}\Delta_i,\varphi_{R*}D_i)$ is a pklt triple. After replacing this triple by $(X_{i+1},\Delta_{i+1},D_{i+1})$, we repeat the process on the triple $(X_{i+1},\Delta_{i+1},D_{i+1})$.

\item If $\varphi_R$ is a $(X_i,\Delta_i,D_i)$-flipping contraction, then by Theorem \ref{thrm:pklt flip exists} there exists a flip $\varphi_i\colon X_i\dashrightarrow X^+_i/X'_i$ of $\varphi_R$. By Proposition \ref{prop:pdisc ineq}, $(X_i^+,\varphi_{i*}\Delta_i,\varphi_{i*}D_i)$ is a pklt triple. After replacing this triple by $(X_{i+1},\Delta_{i+1},D_{i+1})$, we repeat the procedure on the triple $(X_{i+1},\Delta_{i+1},D_{i+1})$.

\item If $\varphi_R$ is not birational, then $(X_i,\Delta_i,D_i)$ with $\varphi_R$ is a Mori fiber space (Definition \ref{def:potential models}).
\end{enumerate}

By repeating this process, we obtain the following sequence of divisorial contractions and flips $\varphi_i$ starting from the triple $(X,\Delta,D)=(X_1,\Delta_1,D_1)$:
\begin{equation}\label{eqn:triple MMP seq}
\begin{tikzcd}
X=X_1 \arrow[dashrightarrow]{r}{\varphi_1} & X_2 \arrow[dashrightarrow]{r}{\varphi_2}&\cdots \arrow[dashrightarrow]{r}{\varphi_{i-1}}&X_i\arrow[dashrightarrow]{r}{\varphi_i}& X_{i+1}\arrow[dashrightarrow]{r}{\varphi_{i+1}}&\cdots%
\end{tikzcd}
\end{equation}
Since each divisorial contraction drops the Picard number by one, any infinite such sequence eventually consists only of flips, i.e., there exists a positive integer $N$ such that $\varphi_i\colon X_i\dashrightarrow X_{i+1}$ is a flip for any $i\geq N$. 

The above results and observations lead to the following theorem.

\begin{theorem}
Let $(X,\Delta,D)$ be a $\Q$-factorial pklt triple. Then we can run the MMP on the triple $(X,\Delta,D)$.
\end{theorem}

However, the following natural question remains to be confirmed in order to guarantee the termination of the MMP on pklt triples.

\begin{conjecture}\label{conj:pklt termination}
Any sequence of flips for a $\Q$-factorial pklt triple $(X,\Delta,D)$ terminates.
\end{conjecture}

We note that this termination does not directly follow from the termination of klt flips as explained in the following remark.

\begin{remark}\label{rem:insufficient pair termination}
Suppose that there exists a sequence of flips for pklt triples:
$$
\begin{tikzcd}
(X_{1},\Delta_{1},D_{1})\arrow[dashrightarrow]{r}{\varphi_{1}}&\cdots\arrow[dashrightarrow]{r}{\varphi_{i-1}}&(X_i,\Delta_i,D_i)\arrow[dashrightarrow]{r}{\varphi_i}& (X_{i+1},\Delta_{i+1},D_{i+1})\arrow[dashrightarrow]{r}{\varphi_{i+1}}&\cdots
\end{tikzcd}
$$
where $\Delta_{i+1}=\varphi_{i*}\Delta_i$ and $D_{i+1}=\varphi_{i*}D_i$.
For each pklt triple $(X_i,\Delta_i,D_i)$, there exists an ample divisor $A_i$ such that $D_i+A_i$ is a $\Q$-divisor and the triple $(X_i,\Delta_i,D_i+A_i)$ is also pklt. By Theorem \ref{thrm:pklt to klt}, there also exists an effective divisor $D'_i\sim_{\R}D_i+A_i$ such that $(X_i,\Delta_i+D'_i)$ is a klt pair. Furthermore, by taking  $A_i$ sufficiently small, we may assume that each $\varphi_i$ is a flip for the klt pair $(X_i,\Delta_i+D'_i)$.
Note that $\Delta_{i+1}+D'_{i+1}$ need not coincide with $\varphi_{i*}(\Delta_i+D'_i)$ by construction and the flips  $\varphi_i$ may depend on essentially distinct klt pairs.
\end{remark}

However, we can prove the following.

\begin{theorem}\label{thrm:klt termination to pklt}
Assume that any sequence of flips of a klt pair terminates. Let $(X,\Delta,D)$ be a $\Q$-factorial pklt triple with a big $\Q$-divisor $D$. Then any sequence of flips of $(X,\Delta,D)$ terminates. 
\end{theorem}
\begin{proof}
Suppose that we are given a sequence of flips for a $\Q$-factorial pklt triple $(X,\Delta,D)=(X_1,\Delta_1,D_1)$. Consider the following sequence of flips as in (\ref{eqn:triple MMP seq}):
$$
\begin{tikzcd}
(X_{1},\Delta_{1},D_{1})\arrow[dashrightarrow]{r}{\varphi_{1}}&\cdots\arrow[dashrightarrow]{r}{\varphi_{i-1}}&(X_i,\Delta_i,D_i)\arrow[dashrightarrow]{r}{\varphi_i}& (X_{i+1},\Delta_{i+1},D_{i+1})\arrow[dashrightarrow]{r}{\varphi_{i+1}}&\cdots.
\end{tikzcd}
$$
By Theorem \ref{thrm:pklt to klt}, there exists an effective divisor $D'_1\sim_{\R}D$ such that $(X,\Delta+D'_1)$ is a klt pair. 
Furthermore, since $\varphi_{i*}(D_i')\sim_{\R}\varphi_{i*}(D_i)=D_{i+1}$ for $i\geq 1$, we may take the effective divisor $D_{i+1}'=\varphi_{i*}(D_i')\sim_{\R} D_{i+1}$ such that $(X_{i+1},\Delta_{i+1}+D'_{i+1})$ is klt. Therefore, the sequence of flips $\varphi_i$ of $(X_i,\Delta_i,D_i)$ can be seen as a sequence of flips starting from the klt pair $(X,\Delta+D_1')=(X_1,\Delta_1+D_1')$:
$$
\begin{tikzcd}
(X_{1},\Delta_{1}+D_{1}')\arrow[dashrightarrow]{r}{\varphi_{1}}&\cdots\arrow[dashrightarrow]{r}{\varphi_{i-1}}&(X_i,\Delta_i+D_i')\arrow[dashrightarrow]{r}{\varphi_i}& (X_{i+1},\Delta_{i+1}+D_{i+1}')\arrow[dashrightarrow]{r}{\varphi_{i+1}}&\cdots.
\end{tikzcd}
$$
By assumption, the sequence of flips of klt pairs terminates and it follows that the given sequence of flips of pklt triples terminates as well.
\end{proof}

\begin{remark}\label{rem:Hashizume termination}
During the preparation of this paper, Hashizume proved the termination of any sequence of flips of the klt pairs $(X,\Delta)$ with big $\Delta$ in \cite{K26}. This result enables us to remove the assumption on the termination of flips in Theorem \ref{thrm:klt termination to pklt}.
\end{remark}

\vspace{1em}

We also consider the MMP on a pklt triple with scaling of an ample divisor.

Let $(X,\Delta,D)=(X_1,\Delta_1,D_1)$ be a $\Q$-factorial pklt triple and $A_1$ an ample divisor such that $K_{X_1}+\Delta_1+D_1+A_1$ is nef. We may assume that $(X_1,\Delta_1,D_1+A_1)$ is pklt.
Define 
$$
\lambda_1\coloneqq \inf\{t\geq 0\mid K_{X_1}+\Delta_1+D_1+tA_1\;\text{is nef}\}.
$$
Suppose that for $i\geq 1$ we have constructed a $\Q$-factorial pklt triple $(X_i,\Delta_i,D_i)$ and the transform $A_i$ of $A_1$ on $X_i$ with a real number $\lambda_{i-1}$ such that $K_{X_i}+\Delta_i+D_i+\lambda_{i-1}A_i$ is nef. We define $\lambda_0=1$.
Define the threshold
$$
\lambda_i\coloneqq \inf\{t\geq0\mid  K_{X_i}+\Delta_i+D_i+t A_i \;\text{is nef}\}.
$$
If $\lambda_i=0$, then $(X_i,\Delta_i,D_i)$ is a minimal model. Otherwise, we have $\lambda_i>0$, which implies that $K_{X_i}+\Delta_i+D_i$ is not nef. Thus, there exists a $(K_{X_i}+\Delta_i+D_i)$-negative extremal ray $R_i$ such that $(K_{X_i}+\Delta_i+D_i+\lambda_i A_i)\cdot R_i=0$. There exists a contraction associated to $R_i$ by Theorem \ref{thrm:triple contraction}. If this contraction gives a Mori fiber space, we terminate the program. Otherwise, the contraction is birational and either it is a divisorial contraction or there exists a flip associated to the ray $R_i$. In either case, we denote the resulting birational map by $\varphi_i\colon X_i\dashrightarrow X_{i+1}$.
By Proposition \ref{prop:pdisc ineq},  $(X_{i+1},\Delta_{i+1},D_{i+1})$ is again a $\Q$-factorial  pklt triple where $\Delta_{i+1}\coloneqq \varphi_{i*}\Delta_i$, $D_{i+1}\coloneqq \varphi_{i*}D_i$. For $A_{i+1}\coloneqq \varphi_{i*}A_i$, we can check that $K_{X_{i+1}}+\Delta_{i+1}+D_{i+1}+\lambda_i A_{i+1}$ is nef.

\begin{theorem}[cf. \protect{\cite[Theorem 1.2]{CJK}}]\label{thrm:pklt-mmp-scaling}
Let $(X,\Delta,D)$ be a $\Q$-factorial pklt triple. Let $A$ be an ample divisor such that $K_X+\Delta+D+A$ is nef. Then we can run the MMP for $(X,\Delta,D)$ with scaling of $A$. Furthermore, if $\Delta+D$ is a big $\Q$-divisor, then every such MMP with scaling terminates.
\end{theorem}
\begin{proof}
By the results we obtained above or \cite[Theorem 1.2]{CJK}, it is easy to see that we can run the MMP on $(X,\Delta,D)$ with scaling of an ample divisor. Thus, it suffices to prove the assertion regarding the termination. 

Suppose that the following is the sequence of steps of the MMP on $(X,\Delta,D)=(X_1,\Delta_1,D_1)$ with scaling of $A$:
$$
\begin{tikzcd}
(X_{1},\Delta_{1},D_{1})\arrow[dashrightarrow]{r}{\varphi_{1}}&\cdots\arrow[dashrightarrow]{r}{\varphi_{i-1}}&(X_i,\Delta_i,D_i)\arrow[dashrightarrow]{r}{\varphi_i}& (X_{i+1},\Delta_{i+1},D_{i+1})\arrow[dashrightarrow]{r}{\varphi_{i+1}}&\cdots.
\end{tikzcd}
$$
Assume that $\Delta+D$ is a big $\Q$-divisor. For every prime divisor $E$ over $X$, we have
    \begin{align*}
        a(E;X,0,\Delta+D)&\geq a(E;X,0)-\ord_E \Delta-\sigma_E(D)\\
        &=a(E;X,\Delta,D).
    \end{align*}
    Thus, $(X,0,\Delta+D)$ is pklt. By Theorem \ref{thrm:pklt to klt}, there is an effective big $\Q$-divisor $B\sim_{\Q} \Delta+D$ such that $(X,B)$ is klt. Every step of the MMP on the triple $(X,\Delta,D)$ is also a step of the MMP on the pair $(X,B)$ since $K_X+B\sim_{\Q} K_X+\Delta+D$ and the linear equivalence is preserved under such steps. Let $A_i$ and $B_i$ be the transforms of $A$ and $B$ on $X_i$. At every step, the nef thresholds of $K_{X_i}+B_i$ and $K_{X_i}+\Delta_i+D_i$ with respect to $A_i$ coincide. Choose an effective divisor $A'\sim_\R A$ such that $(X,B+A')$ is klt. Then the MMP with scaling terminates by \cite[Corollary 1.4.2]{bchm}.
\end{proof}

\bigskip

\section{Admissible rays}\label{sect:admissible rays}

In this section, we study the geometry of extremal rays of the Mori cone $\overline{\NE}(X)$ with respect to pairs and triples. This will provide the key results for the construction of $-(K_X+\Delta)$-MMP.
Unlike the standard MMP, $-(K_X+\Delta)$-MMP requires contracting the $-(K_X+\Delta)$-negative (equivalently, $(K_X+\Delta)$-positive) extremal rays of $\overline{\NE}(X)$. In order to do so, we have to first distinguish which of the $(K_X+\Delta)$-positive extremal rays are even `contractible' and we also have to determine whether such contractions can be regarded as steps of $-(K_X+\Delta)$-MMP. Such criteria will be given in the framework of the MMP on triples. We first define the following extremal rays.

\begin{definition}\label{def:admissible rays}
Let $(X,\Delta)$ be a klt pair. An extremal ray  $R$ of $\overline{\NE}(X)$ is \emph{admissible for} $(X,\Delta)$ if there exists a pklt triple $(X,\Delta,D)$ with a big $\R$-Cartier divisor $D$ such that $R$ is $(K_X+\Delta+D)$-negative. In this case, we say that $D$ \textit{witnesses} the admissible ray $R$ for $(X,\Delta)$.
\end{definition}

Note that for a klt pair $(X,\Delta)$, any $(K_X+\Delta)$-negative extremal ray is admissible for $(X,\Delta)$. (Take a pklt triple $(X,\Delta,D)$ with a sufficiently small ample divisor $D$.) A $(K_X+\Delta)$-positive extremal ray is either  admissible or non-admissible. Therefore, we have
\begin{align*}
\overline{\NE}(X)&=\overline{\sum R_r}\;\;\;\text{($R_r$ are $(K_X+\Delta)$-positive and admissible rays)}\\
&+\overline{\sum R_s}\;\;\;\text{($R_s$ are $(K_X+\Delta)$-positive and non-admissible rays)}\\
&+\overline{\NE}(X)\cap\{\eta\mid \eta\cdot(K_X+\Delta)= 0\}\\
&+\sum R_t\;\;\;\text{($R_t$ are $(K_X+\Delta)$-negative rays)}
\end{align*}
where the closure in the last summation is not necessary due to the cone theorem for klt pairs. See Example \ref{ex:only positive nonadmissible} and Example \ref{ex:-K-MMP on nonMDS} for the examples where all the $(K_X+\Delta)$-positive extremal rays of $\overline{\NE}(X)$ are non-admissible and admissible, respectively. If $R$ is an admissible ray for $(X,\Delta)$, then there exists an associated contraction morphism  $\varphi_R\colon X\to Y$ given by the contraction theorem for pklt triples (Theorem \ref{thrm:triple contraction}) regardless of whether $R$ is $(K_X+\Delta)$-negative or positive.

\medskip

Note that the admissible extremal rays of the Mori cone $\overline{\NE}(X)$ are isolated among the extremal rays by Theorem \ref{thrm:triple cone}. This property for the $(K_X+\Delta)$-negative extremal rays of $\overline{\NE}(X)$ for a klt pair $(X,\Delta)$ is well known. Theorem \ref{thrm:triple cone} also tells us that the extremal rays that are accumulation rays of the extremal rays are not admissible. However, an isolated $(K_X+\Delta)$-positive extremal ray can also be non-admissible (Example \ref{ex:only positive nonadmissible}). 

Recall that if $-(K_X+\Delta)$ is pseudoeffective and $(X,\Delta,-(K_X+\Delta))$ is a pklt triple, we call  $(X,\Delta)$ a pklt pair.

\begin{theorem}\label{thrm:pklt all admissible}
Let $(X,\Delta)$ be a $\Q$-factorial pklt pair with a $\Q$-divisor $\Delta$. Then any $(K_X+\Delta)$-positive extremal ray $R$ is admissible for $(X,\Delta)$ and there exists a rational curve $C$ on $X$ such that $R=\R_{\geq0}[C]$.
\end{theorem}
\begin{proof}
Let $R$ be a $(K_X+\Delta)$-positive extremal ray of $\overline{\NE}(X)$. For any $\gamma\in R\setminus\{0\}$, we have $c\coloneqq \gamma\cdot(K_X+\Delta)>0$. By Theorem \ref{thrm:pklt open} (2), for any sufficiently small $\varepsilon>0$, $(X,\Delta,-(1+\varepsilon)(K_X+\Delta))$ is a pklt triple. Since adding an ample divisor does not increase any $\sigma_E$, we can choose an ample divisor $A$ such that $(X,\Delta,D)$ is pklt where $D=-(1+\varepsilon)(K_X+\Delta)+A$. Note that $D$ is big. Since we can take $A$ arbitrarily small, we may assume that $\gamma\cdot(K_X+\Delta+D)=\gamma\cdot(-\varepsilon(K_X+\Delta)+A)=-\varepsilon c+\gamma\cdot A<0$. Therefore, $R$ is a $(K_X+\Delta)$-positive admissible ray for $(X,\Delta)$. By Theorem \ref{thrm:triple cone} (2), there exists a  rational curve $C$ on $X$ such that $R=\mathbb{R}_{\geq 0} [C]$.
\end{proof}

\begin{theorem}\label{thrm:pklt positive cone}
Let $(X,\Delta)$ be a $\Q$-factorial pklt pair with a $\Q$-divisor $\Delta$. Then there exist at most countably many $(K_X+\Delta)$-positive admissible extremal rays $R_i (i\in I)$ for $(X,\Delta)$ such that 
$$
\overline{\NE}(X)=\overline{\NE}(X)_{K_X+\Delta\leq 0}+\sum_{i\in I} R_i
$$
Furthermore, for any ample divisor $A$, there exist finitely many admissible rays $R_{i_j}$ such that
$$
\overline{\NE}(X)=\overline{\NE}(X)_{K_X+\Delta-A\leq 0}+\sum R_{i_j}
$$
\end{theorem}
\begin{proof}
By the given condition, $(X,\Delta,-(K_X+\Delta))$ is a pklt triple. Thus, by Theorem \ref{thrm:pklt open} (2), the triple $(X,\Delta,-(1+\varepsilon)(K_X+\Delta))$ is pklt for any sufficiently small $\varepsilon>0$. Let $D\coloneqq -(1+\varepsilon)(K_X+\Delta)$. Then, for an extremal ray $R$ of $\overline{\NE}(X)$,  $(K_X+\Delta)\cdot R>0$ if and only if $(K_X+\Delta+D)\cdot R<0$.
Therefore, 
$$\overline{\NE}(X)_{K_X+\Delta+D\geq 0}=\overline{\NE}(X)_{K_X+\Delta\leq 0}$$
and by Theorem \ref{thrm:triple cone}, the first part of Theorem \ref{thrm:pklt positive cone} follows.

For an ample divisor $A$ and $\varepsilon>0$, let $H=\varepsilon A$. Then we have
\begin{align*}
K_X+\Delta+D+H&= K_X+\Delta-(1+\varepsilon)(K_X+\Delta)+H\\
&=-\varepsilon(K_X+\Delta)+H\\
&=-\varepsilon (K_X+\Delta-A).
\end{align*}
Therefore, for an extremal ray $R$ of $\overline{\NE}(X)$, $(K_X+\Delta+D+H)\cdot R<0$ if and only if $(K_X+\Delta-A)\cdot R>0$. Therefore, 
$$
\overline{\NE}(X)_{K_X+\Delta+D+H\geq 0}=\overline{\NE}(X)_{K_X+\Delta-A\leq 0}
$$
and by Theorem \ref{thrm:triple cone}, we obtain the second part of Theorem \ref{thrm:pklt positive cone}.
\end{proof}

Following the convention, we denote 
$$\overline{\NE}(X)_{D_1\geq 0\geq D_2}\coloneqq \overline{\NE}(X)\cap\{\eta\in N_1(X)_\R \mid \eta\cdot D_1\geq 0\geq \eta\cdot D_2\}.$$

\begin{corollary}\label{cor:cone thrm for pklr pair}
Let $(X,\Delta)$ be a $\Q$-factorial pklt pair with a $\Q$-divisor $\Delta$. Then there exist at most countably many $(K_X+\Delta)$-positive admissible rays $R^+_i$ and at most  countably many $(K_X+\Delta)$-negative admissible rays $R^-_j$ such that 
\begin{alignat*}{3}
\overline{\NE}(X)&=\overline{\NE}(X)_{K_X+\Delta\geq0}\;&+\;\overline{\NE}(X)_{K_X+\Delta=0}&+\overline{\NE}(X)_{K_X+\Delta\leq 0}\\
&=\;\;\;\;\;\;\sum R^+_i&+\;\overline{\NE}(X)_{K_X+\Delta=0}&+\;\;\;\;\;\;\sum R^-_j.
\end{alignat*}
Furthermore, for any ample divisor $A$, there exist finitely many $(K_X+\Delta)$-positive admissible rays $R^+_{i_s}$ and finitely many $(K_X+\Delta)$-negative admissible rays $R^-_{j_t}$ such that
$$
\overline{\NE}(X)=\sum R^+_{i_s}+\overline{\NE}(X)_{K_X+\Delta+A\geq 0\geq K_X+\Delta-A}+\sum R^-_{j_t}.
$$
\end{corollary}
\begin{proof}
    This follows from Theorem \ref{thrm:pklt positive cone} and the cone theorem for klt pairs
\end{proof}

Corollary \ref{cor:cone thrm for pklr pair} together with the cone theorem for klt pairs implies that if $(X,\Delta)$ is a $\Q$-factorial pklt pair with a $\Q$-divisor $\Delta$, then the $(K_X+\Delta)$-positive and negative rays can only accumulate toward the hyperplane $\{\gamma\in N_1(X)_{\R}\mid (K_X+\Delta)\cdot\gamma=0\}$. Furthermore, $(K_X+\Delta)$-nontrivial rays are discrete.

\bigskip

Suppose that we are given a $\Q$-factorial klt pair $(X,\Delta)$, where $-(K_X+\Delta)$ is not necessarily pseudoeffective. We discuss how to determine whether a given $(K_X+\Delta)$-positive extremal ray $R$ is admissible for $(X,\Delta)$ or not. We also show that in case $R$ is admissible, the induced birational maps are admissible for $(X,\Delta)$ (Propositions \ref{prop:admissible contraction} and \ref{prop:non-admissible contraction}).

Suppose that there exists a big divisor $D$ such that $D\cdot R<0$. This is a necessary condition for $R$ to be admissible for $(X,\Delta)$.
If no such big divisors exist, then $R$ is not admissible for $(X,\Delta)$ by definition. Define a positive number
$$
\mu_D\coloneqq \sup\{t>0\mid (K_X+\Delta+tD)\cdot R>0\}.
$$
It is clear that $0<\mu_D<+\infty$ and $(K_X+\Delta+\mu_DD)\cdot R=0$.
We also define the threshold (called the plc threshold)
$$
\lambda_D=\textrm{plct}(X,\Delta,D)\coloneqq \sup\{t\geq0\mid (X,\Delta,tD)\;\text{ is pklt}\}.
$$
Now we compare these two thresholds $\mu_D$ and $\lambda_D$. 

\bigskip

\noindent \textbf{Case 1:} If the inequality $\mu_D<\lambda_D$ holds, then for any fixed real number $s$ satisfying $\mu_D< s <\lambda_D$, $(X,\Delta,s D)$ is pklt by Lemma \ref{lem:pklt<plct} and $R$ is $(K_X+\Delta+sD)$-negative. Therefore, $R$ is admissible for $(X,\Delta)$. 

\medskip

Suppose that $R$ is admissible for $(X,\Delta)$ and $R$ is $(K_X+\Delta+D)$-negative for a pklt triple $(X,\Delta,D)$ with a big divisor $D$.
Then the contraction $\varphi_R\colon X\to Z$ for $R$ given by the contraction theorem (Theorem \ref{thrm:triple contraction}) induces a divisorial contraction or a flip, which we denote by $f\colon X\dashrightarrow Y$, with $Y=Z$ in the divisorial case. Note that $\varphi_R$ cannot be a Mori fiber space since $R$ is $D$-negative and $D$ is a big divisor. We prove that the map $f$ is admissible for $(X,\Delta)$ in the sense that the resulting pair  $(Y,\Delta_Y\coloneqq f_{\ast}\Delta)$ is klt.

\begin{proposition}\label{prop:admissible contraction}
Let $(X,\Delta)$, $(X,\Delta,D)$, and $f\colon X\dashrightarrow Y$ be as above. Then $(Y,\Delta_Y)$ is klt.
\end{proposition}
\begin{proof}
Let $\Delta_{Y}, D_{Y}$ be the birational transforms of $\Delta, D$ on $Y$, respectively. By Proposition \ref{prop:pdisc ineq}, we have
$$
a(E;X,\Delta,D)\leq a(E;Y,\Delta_Y,D_Y)
$$
for any prime divisor $E$ over $X$. Since $(X,\Delta,D)$ is pklt, so is $(Y,\Delta_Y,D_Y)$. Furthermore, since
$a(E;Y,\Delta_Y,D_Y)\leq a(E;Y,\Delta_Y)$ holds for any prime divisor $E$ over $X$, $(Y,\Delta_Y)$ is klt.
\end{proof}

\medskip

\noindent\textbf{Case 2:} If $\mu_D<\lambda_D$ does not hold for any big divisor $D$ such that $D\cdot R<0$, then $R$ is not admissible for $(X,\Delta)$. Nonetheless, there can exist a klt pair $(X,\Delta')$ such that $R$ is $(K_X+\Delta')$-negative and there exists a contraction morphism $g=g_R\colon X\to Y$ associated to $R$ by the contraction theorem for the klt pair $(X,\Delta')$. This induces a legitimate step of the MMP for the klt pair $(X,\Delta')$ and the resulting pair is klt when such a step is birational. However, the following shows that $g$ is not admissible for  $(X,\Delta)$.

\begin{proposition}\label{prop:non-admissible contraction}
Let $(X,\Delta)$, $R$, and $g=g_R\colon X\to Y$ be as above.
\begin{enumerate}
\item If $g$ is a divisorial contraction, then $(Y,\Delta_Y)$ is not klt.

\item If $g$ is a small contraction and $f\colon X\dashrightarrow X^+$ is the flip of $g$ with $g^+\colon X^+\to Y$, then $(X^+,\Delta^+\coloneqq f_*\Delta)$ is not klt.

\item If $g$ is a Mori fiber space for $(X,\Delta')$, then $g$ is not a Mori fiber space for $(X,\Delta)$.

\end{enumerate}
\end{proposition}
\begin{proof}
\noindent (1) Suppose that $(Y,\Delta_Y)$ is klt.
The map $g$ applied to $(X,\Delta)$ gives us
$$
K_X+\Delta=g^*(K_Y+\Delta_Y)-aE\;\;\;\text{for some}\;\;\; a>0.
$$
It is easy to see that $(X,\Delta+aE)$ is klt and so is $(X,\Delta+tE)$ for any $t>a$ sufficiently close to $a$. Now choose a sufficiently ample divisor $A$ on $Y$ and let $D\coloneqq g^*A+tE$ which is big. We claim that $(X,\Delta,D)$ is pklt. Let $F$ be a prime divisor over $X$. By \cite[Lemma III.5.14]{Nak}, $\sigma_{F}(D)=\sigma_{F}(A)+t\ord_{F}E=t\ord_{F}E$. Therefore, we have
\begin{align*}
a(F;X,\Delta,D)&=a(F;X,\Delta)-t\ord_{F}E\\
&=a(F;X,\Delta+tE),
\end{align*}
which shows that $(X,\Delta,D)$ is pklt.

Let $C$ be a curve such that $R=\R_{\geq 0}[C]$. Then 
$$
(K_X+\Delta+D)\cdot C=(g^*(K_Y+\Delta_Y+A)+(t-a)E)\cdot C=(t-a)E\cdot C<0
$$
and thus, $R$ is $(K_X+\Delta+D)$-negative. This contradicts the non-admissibility of the ray $R$.

\medskip

\noindent (2) Suppose that $(X^+,\Delta^+)$ is klt and we show that this leads to a contradiction. Take a common resolution $p\colon W\to X$ and $q\colon W\to X^+$ so that we have $h\coloneqq g\circ p=g^+\circ q$. Then by the negativity lemma, we have
$$
p^*(K_X+\Delta)=q^*(K_{X^+}+\Delta^+)-G
$$
for some effective divisor $G$ on $W$ which is both $p,q$-exceptional. This implies
\begin{equation}\label{eqn:pullback flip}
a(F;X,\Delta)=a(F;X^+,\Delta^+)+\ord_F(G)
\end{equation}
for any prime divisor $F$ on $W$. By Lemma \ref{lem:equiv to sim}, $-(K_{X^+}+\Delta^{+})$ is $g^+$-ample. Now choose a sufficiently ample divisor $H$ on $Y$ such that $D^+\coloneqq -(K_{X^+}+\Delta^+)+(g^{+})^*H$ is ample on $X^+$. For a sufficiently small $\varepsilon>0$, let
\begin{align*}
D_{\varepsilon}&\coloneqq -(1+\varepsilon)(K_X+\Delta)+g^*H,\\
D^+_{\varepsilon}&\coloneqq -(1+\varepsilon)(K_{X^+}+\Delta^+)+g^{+*}H.
\end{align*}
Then $D^+_{\varepsilon}$ is ample on $X^+$ by Lemma \ref{lem:equiv to sim} (2) and $D_{\varepsilon}$ is big on $X$.
Since $p^*D_\varepsilon=q^*D^+_{\varepsilon}+(1+\varepsilon)G$, we have
\begin{equation}\label{eqn:pullback flip-2}
\sigma_F(D_{\varepsilon})=(1+\varepsilon)\ord_F(G)
\end{equation}
for every prime divisor $F$ on $W$.
Therefore, by (\ref{eqn:pullback flip}) and (\ref{eqn:pullback flip-2})
\begin{align*}
a(F;X,\Delta,D_{\varepsilon})&=a(F;X,\Delta)-\sigma_F(D_{\varepsilon}) \\
&=a(F;X^+,\Delta^+)+\ord_F(G)-(1+\varepsilon)\ord_F(G)\\
&=a(F;X^+,\Delta^+)-\varepsilon \ord_F(G).
\end{align*}
We can also check that by pulling back to a higher common resolution, the same identities hold for every prime divisor $F$ over $W$, hence over $X$.

Using this equality, we claim that $(X,\Delta,D_{\varepsilon})$ is a pklt triple. Consider a common log resolution $W$ for $(X,\Delta)$ and $(X^+,\Delta^+)$ so that  $K_W+B_W=q^*(K_{X^+}+\Delta^+)$ and $\Supp B_W\cup \Supp G$ is a simple normal crossings divisor. Since $(X^+,\Delta^+)$ is klt, all the coefficients of $B_W$ are $<1$. By taking $\varepsilon>0$ sufficiently small if necessary, the same is true for $B_W+\varepsilon G$. Thus, the pair $(W,B_W+\varepsilon G)$ is sub-klt and 
$$
\inf_F\{a(F;X^+,\Delta^+)-\varepsilon \ord_F(G)\}>0.
$$
Thus, $(X,\Delta, D_{\varepsilon})$ is a pklt triple.
Finally, we can check that
$$
(K_X+\Delta+D_\varepsilon)\cdot R=-\varepsilon (K_X+\Delta)\cdot R<0
$$
which contradicts our assumption that $R$ is not admissible.

\medskip
\noindent (3) The Mori fiber space $g$ is associated to the $(K_X+\Delta')$-negative extremal ray $R$. However, $R$ is $(K_X+\Delta)$-positive and $g$ cannot be a Mori fiber space for $(X,\Delta)$.
\end{proof}

\bigskip

For the admissible rays, much more is true. The following proves Theorem \ref{thrm:elementary -K-MMP}.

\begin{theorem}\label{thrm:-KMMP preserv potdisc}
Let $(X,\Delta)$ be a $\Q$-factorial pair such that $-(K_X+\Delta)$ is pseudoeffective. Let $R$ be a $(K_X+\Delta)$-positive extremal ray of $\overline{\NE}(X)$ which is admissible for $(X,\Delta)$ and $f\colon X\dashrightarrow Y$ the associated divisorial contraction or the flip. Then $-(K_Y+\Delta_Y)=-f_*(K_X+\Delta)$ is pseudoeffective and for any prime divisor $E$ over $X$, we have
$$
a(E;X,\Delta,-(K_X+\Delta))=a(E;Y,\Delta_Y,-(K_Y+\Delta_Y)).
$$
In particular, $(Y,\Delta_Y)$ is again a $\Q$-factorial klt pair such that $-(K_Y+\Delta_Y)$ is pseudoeffective.
Furthermore, for a (partial) $-(K_X+\Delta)$-MMP $\varphi\colon X\dashrightarrow X'$, $(X,\Delta)$ is pklt if and only if $(X',\Delta_{X'})$ is pklt.
\end{theorem}
\begin{proof}
Since $R$ is admissible for $(X,\Delta)$, there exists a big divisor $D$ such that $(X,\Delta,D)$ is pklt and $R$ is $(K_X+\Delta+D)$-negative.
We consider $(X,\Delta)$ as a triple $(X,\Delta,L)$ where $L=-(K_X+\Delta)$. Write $\Delta_Y\coloneqq f_{\ast}\Delta$, $D_Y\coloneqq f_{\ast}D$, and $L_Y\coloneqq f_{\ast}L$. It is clear that $Y$ is $\Q$-factorial and $-(K_Y+\Delta_Y)=-f_*(K_X+\Delta)$ is pseudoeffective.

We first claim that $f$ is $-(K_{X}+\Delta)$-negative. If $f$ is a divisorial contraction, then clearly $f$ is $-(K_{X}+\Delta)$-negative since $f$ is obtained as the contraction of the ray $R$. Suppose that $f\colon X\dashrightarrow Y=X^+$ with $\varphi^+\colon X^+\to T$ is the flip for a $(K_X+\Delta+D)$-negative small contraction $\varphi=\varphi_R\colon X\to T$. Such a flip exists by Theorem \ref{thrm:pklt flip exists}. By the definition of flip, $-(K_X+\Delta+D)$ is $\varphi$-ample and $(K_Y+\Delta_Y+D_Y)$ is $\varphi^+$-ample. Since $-L$ is $\varphi$-ample and $\rho(X/T)=1$, there exists some $c>0$ such that 
$$
-(K_X+\Delta+D)\equiv c(-L)\;\text{ over } T.
$$
Applying $f_{\ast}$ to both sides, we have
$$
-(K_Y+\Delta_Y+D_Y)\equiv c(-L_Y)\;\text{ over } T.
$$
This implies that $-(K_Y+\Delta_Y)$ is $\varphi^+$-ample. Thus, $f$ is also an $L$-flip and $f$ is also $L$-negative.

Now, we prove that equality of the log discrepancies. The following argument is taken from the proof of \cite[Proposition 3.11]{CP}.
Fix a common log resolution $W$ of $(X,\Delta)$ and $(Y,\Delta_Y)$

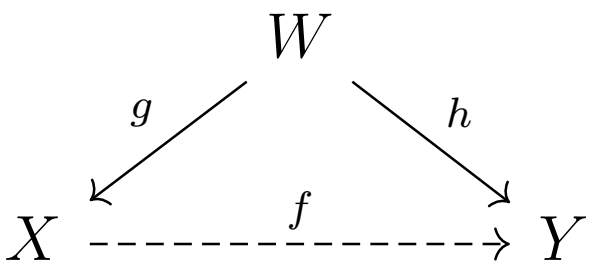
$$\begin{tikzcd}
& W \arrow[dl,"g"'] \arrow[dr,"h"]&\\
X \arrow[rr,dashed,"f"]
&& Y  
\end{tikzcd}
$$
and write
$$
g^*(-(K_X+\Delta))=h^*(-(K_Y+\Delta_Y))+\Gamma
$$
where $\Gamma$ is an $h$-exceptional effective divisor given by
$$
\Gamma=\sum_{E=\text{prime divisor on }W}\big(a(E;X,\Delta)-a(E;Y,\Delta_Y)\big)E.
$$

Let $h^*(-(K_Y+\Delta_Y))=P+N$ be the divisorial Zariski decomposition.
By \cite[Lemma III.5.14]{Nak}, $g^*(-(K_X+\Delta))=P + (N+\Gamma)$ is the divisorial Zariski decomposition.
The following computation shows the required equality:
\begin{align*}
a(E;X,\Delta,-(K_X+\Delta))&=a(E;X,\Delta)-\sigma_E(-(K_X+\Delta))\\
&=a(E;X,\Delta) - \mult_E (N+\Gamma) \\
&=a(E;X,\Delta) - \mult_E N - \big(a(E;X,\Delta)-a(E;Y,\Delta_Y)\big)\\
&=a(E;Y,\Delta_Y) - \sigma_E(-(K_Y+\Delta_Y)) \\
&= a(E;Y,\Delta_Y,-(K_Y+\Delta_Y)). \qedhere
\end{align*}
\end{proof}

\medskip

\section{Minimal Model Program for \texorpdfstring{$-(K_X+\Delta)$}{-(KX + Delta)}}
\label{sect:klt -K MMP}

As an application of the minimal model program on potential triples, we establish the algorithm for the $-(K_X+\Delta)$-MMP on a $\Q$-factorial klt pair $(X,\Delta)$ such that $-(K_X+\Delta)$ is pseudoeffective. Our strategy is to contract $-(K_X+\Delta)$-negative (equivalently, $(K_X+\Delta)$-positive) extremal rays of the Mori cone $\overline{\NE}(X)$. Although the fundamental results for the standard MMP are not applicable here, the cone theorem (Theorem \ref{thrm:triple cone}) and contraction theorem (Theorem \ref{thrm:triple contraction}) for pklt triples come into play.

A naive expectation is that we obtain the following models as the resulting models of the $-(K_X+\Delta)$-MMP.

\begin{definition}\label{def:antimin model}
Let $(X,\Delta)$ be a $\Q$-factorial klt pair such that $-(K_X+\Delta)$ is pseudoeffective.
A birational contraction $\varphi\colon X\dashrightarrow Y$ (or the pair $(Y,\Delta_Y\coloneqq \varphi_*\Delta)$) is called an \emph{anticanonical minimal model of} $(X,\Delta)$ if $(Y,\Delta_Y)$ is a $\Q$-factorial klt pair, $-(K_Y+\Delta_Y)$ is nef, and $\varphi$ is $-(K_X+\Delta)$-nonpositive.
If, additionally, $-(K_Y+\Delta_Y)$ is semiample, then it is called a \emph{good anticanonical minimal model of $(X,\Delta)$}.
\end{definition}

An anticanonical minimal model exists only in special cases, and even when it does exist, one cannot generally expect it to satisfy the goodness property.  For example, if $X$ is the blow-up of $\mathbb P^2$ at nine very general points, then the anticanonical divisor $-K_X$ is nef but not semiample (\cite[Example 6.6]{TX}).
There is also another type of resulting model which we may obtain by running the $-(K_X+\Delta)$-MMP.

\begin{definition}\label{def:nonadmissible model}
Let $(X,\Delta)$ be a $\Q$-factorial klt pair such that $-(K_X+\Delta)$ is pseudoeffective. A birational contraction $\varphi\colon X\dashrightarrow Y$ (or the pair $(Y,\Delta_Y\coloneqq \varphi_*\Delta)$) is called a \emph{non-admissible model} of $(X,\Delta)$ if $\varphi$ is $-(K_{X}+\Delta)$-nonpositive, $(Y,\Delta_Y)$ is a $\Q$-factorial klt pair, $-(K_Y+\Delta_Y)$ is not nef and all the $(K_Y+\Delta_Y)$-positive extremal rays are non-admissible for $(Y,\Delta_Y)$.
\end{definition}

The following gives an example of a non-admissible model.

\begin{example}\label{ex:only positive nonadmissible}
Let $X=\mathbb F_3$ be the Hirzebruch surface of degree $3$, and let $S$ be its negative section with $S^2=-3$. Then $X=\mathbb P_{\mathbb P^1}(E)$ where $E=\mc O_{\mathbb P^1}\oplus \mc{O}_{\mathbb P^1}(3)$.
Take three distinct fibers $F_i (i=1,2,3)$ and let $\Delta\coloneqq a(F_1+F_2+F_3)$ where $\frac{2}{3}<a<1$. Since $K_X=-2S-5F$ for a fiber $F$, we have $-(K_X+\Delta)\sim_{\R} 2S+(5-3a)F$, which is pseudoeffective. Note that $(X,\Delta)$ is klt and $-(K_X+\Delta)$ is not nef since $-(K_X+\Delta)\cdot S=(2S+(5-3a)F)\cdot S=-(1+3a)<0$. One can also check that $R=\R_{\geq 0}[S]$ is the only $(K_X+\Delta)$-positive extremal ray of $\overline{\NE}(X)$.

We can prove that $R$ is a non-admissible ray for $(X,\Delta)$ as follows. Suppose that there exists a big divisor $D$ such that the extremal ray $R$ is $(K_X+\Delta+D)$-negative. Then we claim that $(X,\Delta,D)$ is not pklt. 

Let $D\equiv xS+yF$. Since $D\cdot S=y-3x<0$, the negative part of the Zariski decomposition of $D$ is 
$$N_{\sigma}(D) = \left(x-\frac{y}{3}\right)S = -\frac{D\cdot S}{3}S.$$ Therefore, $\sigma_S(D) = -\frac{D\cdot S}{3} > \frac{1+3a}{3}$. Since $a>\frac{2}{3}$, we have $ \frac{1+3a}{3}>1$. Our assumption $S\not\subset \Supp\Delta$ implies $ a(S;X,\Delta)=1$. Consequently, 
$$
a(S;X,\Delta,D) = a(S;X,\Delta)-\sigma_S(D) = 1-\sigma_S(D) < 0. 
$$
Therefore, the unique $(K_X+\Delta)$-positive extremal ray $R=\mathbb R_{\geq 0}[S]$ cannot be $(K_X+\Delta+D)$-negative for any pklt triple $(X,\Delta,D)$ with a big divisor $D$.

However, there exists a klt pair $(X,\Delta')$ such that $R$ is the unique $(K_X+\Delta')$-negative extremal ray. 
Choose three general curves $C_1,C_2,C_3 \in |S+3F|$. Since $(S+3F)\cdot S=0$, the curves $C_i$ can be chosen disjoint from $S$. One can check that the pair $(X,\Delta')$ is klt where $\Delta' = \frac{1}{2}S + \frac{1}{2}(C_1+C_2+C_3)$. 
Since 
$$
\Delta' \equiv \frac{1}{2}S + \frac{3}{2}(S+3F) = 2S+\frac{9}{2}F,
$$ 
we have
$$
K_X+\Delta' \equiv (-2S-5F) + \left(2S+\frac{9}{2}F\right) = -\frac{1}{2}F.
$$
Therefore, $(K_X+\Delta')\cdot S = -\frac{1}{2}<0$ and $(K_X+\Delta')\cdot F = 0$. This implies that the ray $R=\R_{\geq 0}[S]$ is the unique $(K_X+\Delta')$-negative extremal ray of $\overline{\NE}(X)  = \mathbb{R}_{\geq 0}[S] + \mathbb{R}_{\geq 0}[F]$. By the contraction theorem for klt pairs, the ray $R$ is still contractible by the divisorial contraction $f=f_R\colon X\to Y$ associated to $R$, which is a step of $(K_X+\Delta')$-MMP.

Thus, the pair $(Y,f_{R*}\Delta')$ is klt and one can also check that $Y=\mathbb P(1,1,3)$.
Such contraction $f_R$ of the ray $R$ is not admissible for the given pair $(X,\Delta)$ in the sense that $(Y,f_{R*}\Delta)$ is not klt (actually not even lc). This example illustrates Proposition \ref{prop:non-admissible contraction}. $\hfill\qed$
\end{example}

This gives an example of a variety $X$ having a rational polyhedral $(K_X+\Delta)$-positive half cone $\overline{\NE}(X)_{K_X+\Delta>0}$ and also shows that an isolated and  rational $(K_X+\Delta)$-positive extremal ray can be non-admissible.

\bigskip

Now we are ready to run $-(K_X+\Delta)$-MMP. The following algorithm is called the $-(K_X+\Delta)$-MMP.

Suppose that as an initial model, we are given a $\Q$-factorial klt pair $(X,\Delta)$ such that $-(K_X+\Delta)$ is pseudoeffective. We first ask
\medskip
\begin{center}
\boxed{\textrm{Q1}}: Is $-(K_X+\Delta)$ nef? 
\end{center}
\medskip

If the answer is \emph{yes}, then $(X,\Delta)$ is a $-(K_X+\Delta)$-minimal model and the program terminates.  If the answer is \emph{no}, then $-(K_X+\Delta)$ is not nef and there exist $(K_X+\Delta)$-positive extremal rays of $\overline{\NE}(X)$. We ask
\medskip
\begin{center}
\boxed{\textrm{Q2}}: Among the $(K_X+\Delta)$-positive extremal rays of $\overline{\NE}(X)$, is there an admissible ray?
\end{center}
\medskip
Note that non-isolated $(K_X+\Delta)$-positive extremal rays are not admissible by Theorem \ref{thrm:triple cone}. Thus, in this step, we consider isolated $(K_X+\Delta)$-positive extremal rays if such rays exist. The admissibility of such rays for $(X,\Delta)$ is determined by the threshold comparison explained in Section \ref{sect:admissible rays}. If all the isolated $(K_X+\Delta)$-positive extremal rays are also non-admissible and the answer to $\boxed{\textrm{Q2}}$ is \emph{no}, then $(X,\Delta)$ is a non-admissible model and the $-(K_X+\Delta)$-MMP terminates with this model.

If the answer to $\boxed{\textrm{Q2}}$ is \emph{yes}, then there exists a $(K_X+\Delta)$-positive admissible extremal ray  for $(X,\Delta)$. Then by definition, there exists a big divisor $D$ such that $(X,\Delta,D)$ is pklt and $R$ is $(K_X+\Delta+D)$-negative. By the contraction theorem (Theorem \ref{thrm:triple contraction}) and Theorem \ref{thrm:pklt flip exists}, there exists a divisorial contraction or a flip for the admissible ray $R$ of $(X,\Delta)$, which we denote by $f\colon X\dashrightarrow Y$ in either case.
By Theorem \ref{thrm:-KMMP preserv potdisc}, $(Y,\Delta_Y)$ is klt and $-(K_Y+\Delta_Y)$ is pseudoeffective.  Replace the pair $(X,\Delta)$ by $(Y,\Delta_Y)$ where $\Delta_Y\coloneqq f_*\Delta$ and repeat the process. Note that the process continues to run only if there exists a $(K_X+\Delta)$-positive admissible extremal ray. Furthermore, since each divisorial contraction drops the Picard number $\rho$ by one, we may assume that the MMP on the given klt pair $(X,\Delta)$ consists only of flips $f_i$ of pklt triple $(X_i,\Delta_i,D_i)$:
\begin{equation*}
\begin{tikzcd}
X=X_1 \arrow[dashrightarrow]{r}{f_1} & X_2 \arrow[dashrightarrow]{r}{f_2}&\cdots \arrow[dashrightarrow]{r}{f_{i-1}}&X_i\arrow[dashrightarrow]{r}{f_i}& X_{i+1}\arrow[dashrightarrow]{r}{f_{i+1}}&\cdots%
\end{tikzcd}
\end{equation*}
where $\Delta_i=f_{i-1 *}\Delta_{i-1}$ for $i\geq 2$ and $D_i$ are big divisors witnessing the chosen admissible rays used to construct $f_i$.

We expect that this process terminates after finitely many steps and as a resulting model, we obtain either an anticanonical minimal model or a non-admissible model of $(X,\Delta)$. We remark again that for a non-admissible model $(X',\Delta')$, $-(K_{X'}+\Delta')$ is not nef and Mori fiber spaces cannot occur in our setting.

\begin{conjecture}\label{conj:-KMMP termination}
Let $(X,\Delta)$ be a $\Q$-factorial klt pair such that $-(K_X+\Delta)$ is pseudoeffective. Then the $-(K_X+\Delta)$-MMP terminates after finitely many steps. Equivalently,  there exists no infinite sequence of pklt triples $(X_i,\Delta_i,D_i)$ and their associated flips $f_i$ for a $(K_{X_i}+\Delta_i)$-positive extremal ray which is $(K_{X_i}+\Delta_i+D_i)$-negative:
\begin{equation}\label{spade}
\begin{tikzcd}
X=X_1 \arrow[dashrightarrow]{r}{f_1} & X_2 \arrow[dashrightarrow]{r}{f_2}&\cdots \arrow[dashrightarrow]{r}{f_{i-1}}&X_i\arrow[dashrightarrow]{r}{f_i}& X_{i+1}\arrow[dashrightarrow]{r}{f_{i+1}}&\cdots%
\end{tikzcd}
\end{equation}
where $\Delta_i=(f_{i-1}\circ f_{i-2}\circ\cdots f_1)_*\Delta$ on $X_i$ and $D_i$ are big divisors on $X_i$.
\end{conjecture}

Note that for each flip $f_i$ of $(X_i,\Delta_i,D_i)$, there exists a small ample divisor $A_i$ such that $D_i+A_i$ is a $\Q$-divisor and $f_i$ is also a flip of the triple $(X_i,\Delta_i,D_i+A_i)$. By Theorem \ref{thrm:pklt to klt}, there exists an effective divisor $D'_i\sim_\Q D_i+A_i$ such that $(X_i,\Delta'_i)$ is klt where $\Delta'_i=\Delta_i+D'_i$. 
Thus, each flip $f_i\colon X_i\dashrightarrow X_{i+1}$ for the pklt triple $(X_i,\Delta_i,D_i)$ can be considered as a flip for the klt pair $(X_i,\Delta'_i)$. Therefore, the sequence of flips $f_i$ of triples $(X_i,\Delta_i,D_i)$ of (\ref{spade}) in Conjecture \ref{conj:-KMMP termination} is equivalent to the sequence of flips for klt pairs $(X_i,\Delta'_i)$.
However, there is no reason to expect $f_{i*}\Delta'_{i-1}=\Delta'_i$ to hold for all $i$.

\medskip

\bigskip

Next, we prove the dichotomy theorem (Theorem \ref{thrm:dichotomy}) which enables us to determine the type of resulting model of the $-(K_X+\Delta)$-MMP.

Suppose that we are given a $\Q$-factorial klt pair $(X,\Delta)$ with a $\Q$-divisor $\Delta$ such that $-(K_X+\Delta)$ is pseudoeffective.

Recall first that if $(X,\Delta)$ is pklt and $-(K_X+\Delta)$ is big, then $X$ is of Fano type by \cite[Corollary 3.10]{CJK}, hence a Mori dream space by \cite{bchm}. Suppose that $f\colon X\dashrightarrow X'$ is the $D$-MMP for $D=-(K_X+\Delta)$ which gives us a $D$-nef model $X'$ and $D_{X'}\coloneqq f_*D=-(K_{X'}+\Delta_{X'})$ is semiample. Since $(X,\Delta)$ is pklt, Theorem \ref{thrm:-KMMP preserv potdisc} implies that $(X',\Delta_{X'})$ is klt. Therefore, $(X',\Delta_{X'})$ is a good anticanonical minimal model of $(X,\Delta)$.
\bigskip

\noindent \textbf{Case 1}: Suppose that the given pair $(X,\Delta)$ is pklt and $-(K_X+\Delta)$ is not nef. Then by Theorem \ref{thrm:pklt all admissible}, all $(K_X+\Delta)$-positive extremal rays are admissible for $(X,\Delta)$. Let $R$ be a $(K_X+\Delta)$-positive extremal ray and denote its associated elementary birational map by $f\colon X\dashrightarrow Y$. By Theorem \ref{thrm:-KMMP preserv potdisc}, $(Y,\Delta_Y)$ is also pklt and $-(K_Y+\Delta_Y)$ is pseudoeffective. Therefore, as we run $-(K_X+\Delta)$-MMP, we only encounter the admissible rays as we run the $-(K_X+\Delta)$-MMP. If we assume Conjecture \ref{conj:-KMMP termination}, then after finitely many steps we obtain an anticanonical minimal model of $(X,\Delta)$.

\smallskip

\noindent \textbf{Case 2}: Suppose that the given pair $(X,\Delta)$ is not pklt and $-(K_X+\Delta)$ is not nef. In this case, a $(K_X+\Delta)$-positive extremal ray can be either admissible or non-admissible. (By Example \ref{ex:only positive nonadmissible}, there can be only non-admissible rays in $\overline{\NE}(X)_{K_X+\Delta>0}$.) If there exists an admissible $(K_X+\Delta)$-positive extremal ray $R$ for $(X,\Delta)$ and we denote its associated elementary birational map by $f\colon X\dashrightarrow Y$, then $(Y,\Delta_Y)$ is not pklt by Theorem \ref{thrm:-KMMP preserv potdisc} and $-(K_Y+\Delta_Y)$ is pseudoeffective. Thus, we can repeat the same process with $(Y,\Delta_Y)$. Assuming that Conjecture \ref{conj:-KMMP termination} holds, after finitely many steps we arrive at a non-admissible model $(X',\Delta_{X'})$ of $(X,\Delta)$ which has only non-admissible extremal rays in $\overline{\NE}(X')_{K_{X'}+\Delta_{X'}>0}$. 
    
\medskip

This observation leads to the following theorem.

\begin{theorem}\label{thrm:dichotomy}
Let $(X,\Delta)$ be a $\Q$-factorial klt pair with a $\Q$-divisor $\Delta$ such that $-(K_X+\Delta)$ is pseudoeffective.
Assume that Conjecture \ref{conj:-KMMP termination} holds.
\begin{enumerate}[leftmargin=18pt]
    \item The pair $(X,\Delta)$ is pklt if and only if there exists an anticanonical minimal model of $(X,\Delta)$. 
    \item The pair $(X,\Delta)$ is not pklt if and only if there exists a non-admissible model of $(X,\Delta)$.
\end{enumerate}
\end{theorem}
\begin{proof}
The `only if' directions for both cases (1) and (2) are proved above under the assumption that $-(K_X+\Delta)$-MMP terminates (Conjecture \ref{conj:-KMMP termination}).
The `if' directions follow from Theorem \ref{thrm:-KMMP preserv potdisc}. 
\end{proof}

The following gives an example of a variety that is not a Mori dream space and admits a terminating $-(K_X+\Delta)$-MMP. It was studied previously by Totaro in \cite{To10} and we re-examine it in the framework of the $-(K_X+\Delta)$-MMP (see also \cite[Example 4.1]{KL2}).

\begin{example}\label{ex:-K-MMP on nonMDS} (cf. \cite[\S 2, Example]{To10})
Let $X$ be the blowup of $\mathbb P^2$ at the following twelve points:
$$
P_0=[1:0:0], P_1=[0:1:0], P_2=[0:0:1], P_{i,j}=[1:\zeta^i:\zeta^j] \;(i,j=0,1,2)
$$
where $\zeta$ is a primitive cubic root of unity.
(This is the dual of the well-known Hesse configuration.) 
Let $E_s$ ($s=1,2,\cdots,12$) be the exceptional divisors on $X$ over the points $P_{k}$ and $P_{i,j}$.
Let $L_t$ ($t=1,2,\cdots, 9$) be the nine lines in $\mathbb P^2$ through three of $P_{i,j}$ and one of $P_k$.
Let $\overline{L}_t$ be the strict transforms of $L_t$ on $X$. The rays $R_t=\mathbb R_{\geq 0}[\overline{L_t}]$ form $K_X$-positive extremal rays of the Mori cone $\overline{\NE}(X)$.
For example, let $L_t=\{y-z=0\}$ on $\mathbb P^2$ for some $t$, and consider its strict transform $\overline{L_t}$ on $X$.
Then we can check that $\overline{L_t}^2=-3$ and $K_X\cdot \overline{L_t}>0$ (i.e., $\overline{L_t}$ is $K_X$-positive).

Now consider two distinct lines in $\mathbb P^2$ which do not pass through any of the twelve points $P_k$ and $P_{i,j}$.
Let $C_1,C_2$ be their strict transforms in $X$. 
Let $\overline{L}\coloneqq \sum_{t=1}^9\overline{L_t}$
and consider 
$$D=a\overline{L}+\frac{1}{2}C_1+\frac{1}{2}C_2$$ 
for $\frac{2}{3}<a<1$.
Let $H$ be the pullback of a line in $\P^2$. Since $D\equiv H+a\overline{L}$, with $H$ big and $a\overline{L}$ effective, the divisor $D$ is big.
Since $\overline{L_t}$ are pairwise disjoint and we may assume that $C_1,C_2$ and $\overline{L}$ have simple normal crossings support, $(X,D)$ is klt. 
We can check that each $\overline{L_t}$ is $K_X$-positive, and $(K_X+D)$-negative:
$$
(K_X+D)\cdot \overline{L_t}=1-3a+1<0 \;\;\text{for any } t=1,2,\cdots, 9.
$$
We can also check that $D$ has the Zariski decomposition with the negative part $N=(a-\frac{1}{3})\overline{L}$.
This shows that the triple $(X,0,D)$ is pklt. Therefore, each ray $\R_{\geq0}[\overline{L_t}]$ is an admissible $K_X$-positive extremal ray of $\overline{\NE}(X)$ and the divisor $D$ witnesses all the $K_X$-positive admissible rays $\R_{\geq 0}[\overline{L_t}]$ for $(X,0)$. Each of the configuration lines $L_t$ contains four of the twelve points, and each point lies on three of the lines $L_t$. Consequently, we have
$$
\overline{L} \sim 9H - 3\sum_{s=1}^{12} E_s = -3K_X,
$$
which also implies that $-K_X$ is pseudoeffective.

By Theorem \ref{thrm:triple contraction}, we can contract each $\overline{L_t}$ and the resulting variety is klt. Since $\overline{L_t}$ are disjoint, the admissible contraction for the chosen ray $\R_{\geq 0}[\overline{L_t}]$ does not change the local geometry near the remaining $\overline{L_{t'}}$. Therefore, we can contract all $\overline{L_t}$ and obtain a map $X\to Y$ with $3K_Y\sim 0$.
Thus, $-K_Y$ is semiample  and $Y$ is a good anticanonical minimal model of $X$. However, since the nef cone $\Nef(Y)$ is circular in $N^1(Y)_{\R}\cong \mathbb R^4$, the variety $X$ is not a Mori dream space. 
\end{example}

\medskip

\clearpage

Lastly, we summarize the $-(K_X+\Delta)$-MMP by giving the following flow chart.

\bigskip

\clearpage

\bibliographystyle{habbvr}
\bibliography{biblio}

\end{document}